\documentclass[twoside,a4paper,10pt]{article}
\usepackage{latexsym,amssymb}
\usepackage{mltex}
\usepackage{hyperref}
\usepackage{graphicx}
\usepackage[frenchb,english]{babel}
\usepackage[matrix,arrow]{xy}
\def \be{\begin{eqnarray*}}
\def \ee{\end{eqnarray*}}
\def \ben{\begin{enumerate}}
\def \een{\end{enumerate}}
\def \beit{\begin{itemize}}
\def \eeit{\end{itemize}}
\def \bui#1#2{\mathrel{\mathop{\kern 0pt#1}\limits^{#2}}}
\def \buil#1#2{\mathrel{\mathop{\kern 0pt#1}\limits_{#2}}}
\def \bfll{\begin{flushleft}}
\def \efll{\end{flushleft}}
\def \bflr{\begin{flushright}}
\def \eflr{\end{flushright}}

\def \findemo{\hfill$\square$}
\def \lra{\longrightarrow}
\def \lmt{\longmapsto}
\def \ovl{\overline}
\def \wih{\widehat}
\def \wit{\widetilde}
\def \wnabla{\wit{\nabla}}
\def \cdotM{\buil{\cdot}{M}}

\def \C{\mathbb{C}}

\def \R{\mathbb{R}}
\def \tr{\mathrm{tr}}

\def \grad{\mathrm{grad}}
\def \S{\mathbb{S}}

\def \la{\langle}
\def \ra{\rangle}

\newcommand{\pa}[1]{\left(#1\right)}

\newtheorem{ethm}{Theorem}[section]
\newtheorem{edefi}[ethm]{Def\mbox{}inition}
\newtheorem{elemme}[ethm]{Lemma}
\newtheorem{erem}[ethm]{Remark}
\newtheorem{erems}[ethm]{Remarks}
\newtheorem{ecor}[ethm]{Corollary}
\newtheorem{prop}[ethm]{Proposition}

\title{Imaginary K\"ahlerian Killing spinors I}
\author{Nicolas Ginoux\footnote{Fakult\"at f\"ur Mathematik,
Universit\"at Regensburg,
D-93040 Regensburg,
E-mail: \texttt{nicolas.ginoux@mathematik.uni-regensburg.de}}  and Uwe Semmelmann\footnote{Institut f\"ur Geometrie und Topologie,
Universit\"at Stuttgart,
Pfaffenwaldring 57,
D-70569 Stuttgart, E-mail: \texttt{uwe.semmelmann@mathematik.uni-stuttgart.de}}}

\begin{document}
\maketitle

\noindent\begin{center}\begin{tabular}{p{155mm}}
\begin{small}{\bf Abstract.} We describe and to some extent characterize a new family of K\"ahler spin manifolds admitting non-trivial imaginary K\"ahlerian Killing spinors.
\end{small}\\
\end{tabular}\end{center}

\section{Introduction}\label{s:intro}

\noindent Let $(\wit{M}^{2n},g,J)$ a K\"ahler manifold of real dimension $2n$ and with K\"ahler-form $\wit{\Omega}$ defined by $\wit{\Omega}(X,Y):=g(J(X),Y)$
for all vectors $X,Y\in T\wit{M}$.
We denote by $p_+:TM\lra T^{1,0}M$, $X\mapsto\frac{1}{2}(X-iJ(X))$ and $p_-:TM\lra T^{0,1}M$, $X\mapsto\frac{1}{2}(X+iJ(X))$ the projection maps.
In case $\wit{M}^{2n}$ is spin, we denote its complex spinor bundle by $\Sigma\wit{M}$.

\begin{edefi}
Let $(\wit{M}^{2n},g,J)$ a spin K\"ahler manifold and $\alpha\in\mathbb{C}$.
A pair $(\psi,\phi)$ of sections of $\Sigma\wit{M}$ is called an \emph{$\alpha$-K\"ahlerian Killing spinor} if and only if it satisf\mbox{}ies, for every $X\in\Gamma(T\wit{M})$,
\[\left|\begin{array}{ll}\wnabla_X\psi&= -\alpha p_-(X)\cdot\phi\\
\wnabla_X\phi&= -\alpha p_+(X)\cdot\psi.\end{array}\right.\]
An $\alpha$-K\"ahlerian Killing spinor is said to be \emph{real} (resp. \emph{imaginary}) if and only if $\alpha\in\R$ (resp. $\alpha\in i\R^*$).
\end{edefi}

\noindent If $\alpha=0$, then an $\alpha$-K\"ahlerian Killing spinor is nothing but a pair of parallel spinors.
The classif\mbox{}ication of K\"ahler spin manifolds (resp. spin manifolds) admitting real non-parallel K\"ahlerian Killing (resp. parallel) spinors has been established by A. Moroianu in \cite{Moroi95} (resp. by McK. Wang in \cite{Wang89}).\\

\noindent In this paper, we describe and partially classify those K\"ahler spin manifolds carrying non-trivial imaginary K\"ahlerian Killing spinors.
Note first that there is no restriction in assuming $\alpha=i$: obviously, changing $(\psi,\phi)$ into $(\psi,-\phi)$ changes $\alpha$ into $-\alpha$; moreover, $(\psi,\phi)$ is an $\alpha$-K\"ahlerian Killing spinor on $(\wit{M}^{2n},g,J)$ if and only if it is an $\frac{\alpha}{\lambda}$-K\"ahlerian Killing spinor on $(\wit{M}^{2n},\lambda^2g,J)$ for any constant $\lambda>0$.\\

\noindent K.-D. Kirchberg, who introduced this equation (see \cite{KirchbergAGAG93} for references), showed that, if a non-zero $i$-K\"ahlerian Killing spinor $(\psi,\phi)$ exists on $(\wit{M}^{2n},g,J)$, then necessarily the complex dimension $n$ of $\wit{M}$ is odd, the manifold $(\wit{M}^{2n},g)$ is Einstein with scalar curvature $-4n(n+1)$, the pair $(\psi,\phi)$ vanishes nowhere and satisf\mbox{}ies $\wit{\Omega}\cdot\psi=-i\psi$ as well as $\wit{\Omega}\cdot\phi=i\phi$, see \cite{KirchbergAGAG93} and Proposition \ref{pKirch} below for further properties.
Moreover, he proved in the case $n=3$ that the holomorphic sectional curvature must be constant \cite[Thm. 16]{KirchbergAGAG93}, in particular only the complex hyperbolic space $\C\mathrm{H}^3$ occurs as simply-connected complete $(\wit{M}^{6},g,J)$ with non-trivial $i$-K\"ahlerian Killing spinors.\\

\noindent We extend Kirchberg's results in several ways.
First, we study in detail the critical points of the length function $|\psi|$ of $\psi$.
We show that, if the underlying Riemannian manifold $(\wit{M}^{2n},g)$ is connected and complete, then $|\psi|$ has at most one critical value, which then has to be a (global) minimum and that the corresponding set of critical points is a K\"ahler totally geodesic submanifold (Proposition \ref{pcrit}).\\
As a next step, we describe a whole family of examples of K\"ahler manifolds admitting non-trivial $i$-K\"ahlerian Killing spinors (Theorem \ref{pclassifdwp}), including the complex hyperbolic space and some K\"ahler ma\-ni\-folds with non-constant holomorphic sectional curvature (Corollary \ref{c:exnoncstholsect}).
All arise as so-called doubly-warped products over Sasakian manifolds.
A more detailed study of the induced spinor equation on that Sasakian manifold allows the complex hyperbolic space to be characterized within the family (Theorem~\ref{pcarCHn}).\\

\noindent In the last section, we show that doubly-warped products are the only possible K\"ahler manifolds with non-trivial $i$-K\"ahlerian Killing spinors as soon as both components of $(\psi,\phi)$ have the same length and are exchanged through the Clifford multiplication by a (real) vector field (Theorem \ref{pclassiftype1}).
This shows an interesting analogy with H. Baum's classification \cite{Baum891,Baum892} of complete Riemannian spin manifolds with imaginary Killing spinors.\\

\noindent{\bf\small Acknowledgment.} {\small This project benefited from the generous support of the universities of Hamburg, Potsdam, Cologne and Regensburg as well as the DFG-Sonderforschungsbereich 647.
Special thanks are due to Christian B\"ar and Bernd Ammann.
We also acknowledge very helpful discussions with Bogdan Alexandrov, Georges Habib and Daniel Huybrechts.}

\section{General integrability conditions}\label{s:genintegcond}

\noindent In this section we look for further necessary conditions for the existence of imaginary K\"ahlerian Killing spinors.
Consider the vector f\mbox{}ield $V$ on $\wit{M}$ def\mbox{}ined by
\begin{equation}\label{eqdefV}g(V,X):=\Im m(\langle p_+(X)\cdot\psi,\phi\rangle)\end{equation} for every vector $X$ on $\wit{M}$.
We recall the following

\begin{prop}[see \cite{KirchbergAGAG93}]\label{pKirch}
Let $(\psi,\phi)$ be an $i$-K\"ahlerian Killing spinor on $(\wit{M}^{2n},g,J)$ which does not vanish identically.
Then the following properties hold:
\beit\item[i)] $\grad(|\psi|^2)=\grad(|\phi|^2)=2V$.
\item[ii)] For all vectors $X,Y\in T\wit{M}$, 
\[g(\wnabla_XV,Y)=\Re e\pa{\langle p_-(X)\cdot\phi,p_-(Y)\cdot\phi\rangle+\langle p_+(X)\cdot\psi,p_+(Y)\cdot\psi\rangle}.\]
In particular,
\[\mathrm{Hess}(|\psi|^2)(X,Y)=\mathrm{Hess}(|\phi|^2)(X,Y)=2\Re e\pa{\langle p_-(X)\cdot\phi,p_-(Y)\cdot\phi\rangle+\langle p_+(X)\cdot\psi,p_+(Y)\cdot\psi\rangle}.\]
\item[iii)] $\Delta(|\psi|^2)=\Delta(|\phi|^2)=-2(n+1)(|\psi|^2+|\phi|^2)$, where $\Delta:=-\tr_g(\mathrm{Hess})$.
\item[iv)] The vector f\mbox{}ield $V$ is holomorphic, i.e., it satisf\mbox{}ies: $\wnabla_{J(X)}V=J(\wnabla_XV)$ for every $X\in T\wit{M}$. In particular, the vector f\mbox{}ield $J(V)$ is Killing on $\wit{M}$.
\item[v)] $\grad(|V|^2)=2\wnabla_VV$.
\eeit
\end{prop}

\noindent Note that, from Proposition \ref{pKirch}, the identity $\Delta(|\psi|^2+|\phi|^2)=-4(n+1)(|\psi|^2+|\phi|^2)$ holds on $\wit{M}$,
therefore $\wit{M}$ cannot be compact.\\

\noindent Next we are interested in the critical points of $|\psi|^2$ (or of $|\phi|^2$, they are the same by Proposition \ref{pKirch}$.i)$).
We need a technical lemma:

\begin{elemme}\label{lnabla2V}
Under the hypotheses of {\rm Proposition \ref{pKirch}}, one has
\[\wnabla_X\wnabla_YV=\wnabla_{\wnabla_XY}V+\{2g(V,X)Y+g(V,Y)X-g(V,J(Y))J(X)+g(X,Y)V+g(J(X),Y)J(V)\}\]
for all vector f\mbox{}ields $X,Y$ on $\wit{M}$.
Therefore,
\[\mathrm{Hess}(|V|^2)(X,Y)=2g(\wnabla_XV,\wnabla_YV)+2\pa{3g(X,V)g(Y,V)+|V|^2g(X,Y)-g(X,J(V))g(Y,J(V))}.\]
\end{elemme}

\noindent{\it Proof}: Using Proposition \ref{pKirch}, we compute in a local orthonormal basis $\{e_j\}_{1\leq j\leq 2n}$ of $T\wit{M}$:
\be
\wnabla_X\wnabla_YV&=&\sum_{j=1}^{2n}\Re e\Big(\langle p_-(\wnabla_XY)\cdot\phi,p_-(e_j)\cdot\phi\rangle+\langle p_+(\wnabla_XY)\cdot\psi,p_+(e_j)\cdot\psi\rangle\\
& &\phantom{\sum_{j=1}^{2n}\Re e\Big(}+\langle p_-(Y)\cdot\wnabla_X\phi,p_-(e_j)\cdot\phi\rangle+\langle p_-(Y)\cdot\phi,p_-(e_j)\cdot\wnabla_X\phi\rangle\\
& &\phantom{\sum_{j=1}^{2n}\Re e\Big(}+\langle p_+(Y)\cdot\wnabla_X\psi,p_+(e_j)\cdot\psi\rangle+\langle p_+(Y)\cdot\psi,p_+(e_j)\cdot\wnabla_X\psi\rangle\Big)e_j\\
&=&\sum_{j=1}^{2n}\Re e\Big(\langle p_-(\wnabla_XY)\cdot\phi,p_-(e_j)\cdot\phi\rangle+\langle p_+(\wnabla_XY)\cdot\psi,p_+(e_j)\cdot\psi\rangle\\
& &\phantom{\sum_{j=1}^{2n}\Re e\Big(}-\alpha\langle p_-(Y)\cdot p_+(X)\cdot\psi,p_-(e_j)\cdot\phi\rangle+\alpha\langle p_-(Y)\cdot\phi,p_-(e_j)\cdot p_+(X)\cdot\psi\rangle\\
& &\phantom{\sum_{j=1}^{2n}\Re e\Big(}-\alpha\langle p_+(Y)\cdot p_-(X)\cdot\phi,p_+(e_j)\cdot\psi\rangle+\alpha\langle p_+(Y)\cdot\psi,p_+(e_j)\cdot p_-(X)\cdot\phi\rangle\Big)e_j\\
&=&\wnabla_{\wnabla_XY}V\\
& &+\sum_{j=1}^{2n}\Im m\Big(\langle p_-(Y)\cdot p_+(X)\cdot\psi,p_-(e_j)\cdot\phi\rangle+\langle p_+(Y)\cdot p_-(X)\cdot\phi,p_+(e_j)\cdot\psi\rangle\Big)e_j\\
& &-\sum_{j=1}^{2n}\Im m\Big(\langle p_-(Y)\cdot\phi,p_-(e_j)\cdot p_+(X)\cdot\psi\rangle+\langle p_+(Y)\cdot\psi,p_+(e_j)\cdot p_-(X)\cdot\phi\rangle\Big)e_j.
\ee
We compute the second line of the right-hand side of the preceding equation (the treatment of the third one is analogous).
Using $\langle p_+(X)\cdot\psi,\phi\rangle=2ig(V,p_+(X))$, we obtain 
\be
\langle p_+(Y)\cdot p_-(X)\cdot\phi,p_+(e_j)\cdot\psi\rangle&=&\ovl{\langle\psi,p_-(X)\cdot p_+(Y)\cdot p_-(e_j)\cdot\phi\rangle}+4ig(Y,p_-(e_j))g(V,p_-(X))\\
& &+4ig(Y,p_-(X))g(V,p_-(e_j)).
\ee
We deduce that, for every $j\in\{1,\ldots,2n\}$,
\be\langle p_-(Y)\cdot p_+(X)\cdot\psi,p_-(e_j)\cdot\phi\rangle+\langle p_+(Y)\cdot p_-(X)\cdot\phi,p_+(e_j)\cdot\psi\rangle
&=&\hspace{-3mm}2\Re e\pa{\langle\psi,p_-(X)\cdot p_+(Y)\cdot p_-(e_j)\cdot\phi\rangle}\\
& &+4ig(Y,p_-(e_j))g(V,p_-(X))\\
& &+4ig(Y,p_-(X))g(V,p_-(e_j)).
\ee
The imaginary part of the right-hand side of the last equality is then given for every $j\in\{1,\ldots,2n\}$ by
\be
4\Re e\pa{g(Y,p_-(e_j))g(V,p_-(X))+g(Y,p_-(X))g(V,p_-(e_j))}
&=&g(V,X)g(Y,e_j)+g(V,J(X))g(J(Y),e_j)\\
& &\hspace{-4mm}+g(X,Y)g(V,e_j)+g(J(X),Y)g(J(V),e_j).
\ee
This shows that
\be
\sum_{j=1}^{2n}\Im m\Big(\langle p_-(Y)\cdot p_+(X)\cdot\psi,p_-(e_j)\cdot\phi\rangle+\langle p_+(Y)\cdot p_-(X)\cdot\phi,p_+(e_j)\cdot\psi\rangle\Big)e_j&=&g(V,X)Y\\
& &+g(V,J(X))J(Y)\\
& &+g(X,Y)V\\
& &+g(J(X),Y)J(V).
\ee
Similarly, one shows that
\be
\sum_{j=1}^{2n}\Im m\Big(\langle p_-(Y)\cdot\phi,p_-(e_j)\cdot p_+(X)\cdot\psi\rangle+\langle p_+(Y)\cdot\psi,p_+(e_j)\cdot p_-(X)\cdot\phi\rangle\Big)e_j&=&-g(V,Y)X\\
& &+g(V,J(Y))J(X)\\
& &-g(V,X)Y\\
& &+g(V,J(X))J(Y).
\ee
Combining the computations above, we obtain
\be
\wnabla_X\wnabla_YV&=&\wnabla_{\wnabla_XY}V\\
& &+\pa{g(V,X)Y+g(V,J(X))J(Y)+g(X,Y)V+g(J(X),Y)J(V)}\\
& &-\pa{-g(V,Y)X+g(V,J(Y))J(X)-g(V,X)Y+g(V,J(X))J(Y)}\\
&=&\wnabla_{\wnabla_XY}V\\
& &+\pa{2g(V,X)Y+g(V,Y)X-g(V,J(Y))J(X)+g(X,Y)V+g(J(X),Y)J(V)},
\ee
which shows the f\mbox{}irst identity.
We deduce for the Hessian of $|V|^2$ that, for all vector f\mbox{}ields $X,Y$ on $\wit{M}$,
\be
\mathrm{Hess}(|V|^2)(X,Y)&=&2g(\wnabla_X\wnabla_VV,Y)\\
&=&2g(\wnabla_{\wnabla_XV}V,Y)+2\Big(2g(V,X)g(V,Y)+|V|^2g(X,Y)-0+g(X,V)g(V,Y)\\
& &\phantom{2g(\wnabla_{\wnabla_XV}V,Y)+2\Big(}+g(J(X),V)g(J(V),Y)\Big)\\
&=&2g(\wnabla_XV,\wnabla_YV)+2\pa{3g(X,V)g(Y,V)+|V|^2g(X,Y)-g(X,J(V))g(Y,J(V))},
\ee
which is the second identity.
This concludes the proof of Lemma \ref{lnabla2V}.
\findemo
$ $\\

\noindent We can now describe more precisely the set of critical values and points of $|\psi|^2$ and $|V|^2$.



\begin{prop}\label{pcrit}
Under the hypotheses of {\rm Proposition \ref{pKirch}}, assume furthermore $(\wit{M}^{2n},g)$ to be connected and complete.
Then the following holds:
\beit\item[i)] The set $\{V=0\}$ of zeros of $V$ coincides with $\{\wnabla_VV=0\}$.
As a consequence, the zeros of $V$ are the only critical points of the function $|V|^2$ on $\wit{M}^{2n}$. 
\item[ii)] The subset $\{V=0\}$ is a (possibly empty) connected totally geodesic K\"ahler submanifold of complex dimension $k<n$ in $(\wit{M}^{2n},g,J)$.
Furthermore, for all $x,y\in\{V=0\}$, every geodesic segment between $x$ and $y$ lies in $\{V=0\}$.
\item[iii)] The function $|\psi|^2$ has at most one critical value on $\wit{M}^{2n}$, which is then a global minimum of $|\psi|^2$.
Furthermore, the set of critical points of $|\psi|^2$ is a connected totally geodesic K\"ahler submanifold in $(\wit{M}^{2n},g,J)$.
\eeit
\end{prop}

\noindent{\it Proof}: The proof relies on simple computations and arguments.\\
$i)$ Proposition \ref{pKirch}$.v)$ already implies that $\{\wnabla_VV=0\}$ coincides with the set of critical points of $|V|^2$.
Every zero of $V$ is obviously a zero of $\wnabla_VV$, i.e., a critical point of $|V|^2$.
Conversely, let $x\in\{\wnabla_VV=0\}$.
Then $0=g_x(\wnabla_VV,V)=|p_-(V_x)\cdot\phi|^2+|p_+(V_x)\cdot\psi|^2$, so that $p_-(V_x)\cdot\phi=0$ and $p_+(V_x)\cdot\psi=0$, which, in turn, implies $0=\Im m\pa{\langle p_+(V_x)\cdot\psi,\phi\rangle}=g(V_x,V_x)$, that is, $V_x=0$.
This shows $i)$.\\
$ii)$ 
The subset $\{V=0\}$ - if non-empty - is the fixed-point-set in $\wit{M}^{2n}$ of the flow of the holomorphic Killing field $J(V)$, therefore it is a totally geodesic K\"ahler submanifold of $\wit{M}^{2n}$ (see e.g. \cite[Sec. II.5]{Kocomplexgeom}); moreover, it cannot contain any open subset of $\wit{M}^{2n}$ since otherwise $V$ would identically vanish as a holomorphic vector field.
To show the connectedness of $\{V=0\}$, it suffices to prove the second part of the statement.
Pick any two points $x_0,x_1$ in $\{V=0\}$ (or, equivalently, any critical points of $|V|^2$) and any geodesic $c$ in $(\wit{M}^{2n},g)$ with $c(0)=x_0$ and $c(1)=x_1$.
Consider the real-valued function $f(t):=|V|_{c(t)}^2$ def\mbox{}ined on $\R$.
Then, for any $t\in\R$ one has $f'(t)=g(\grad(|V|^2),c'(t))=2g(\wnabla_{c'(t)}V,V)$ and
\[ f''(t)=\mathrm{Hess}(|V|^2)(c'(t),c'(t)).\]
Lemma \ref{lnabla2V} provides the Hessian of $|V|^2$: for every $X\in T\wit{M}$,
\[ \mathrm{Hess}(|V|^2)(X,X)=2|\wnabla_XV|^2+2\pa{3g(V,X)^2+|V|^2|X|^2-g(X,J(V))^2}.\]
By Cauchy-Schwarz inequality, $|V|^2|X|^2-g(X,J(V))^2\geq 0$, so that $\mathrm{Hess}(|V|^2)(X,X)\geq 0$ for all $X$, in particular $f$ is convex.
This in turn implies that, if $f'(0)=f'(1)=0$, then necessarily $f$ vanishes on $[0,1]$. 
This proves $ii)$.\\
$iii)$ Set, for any $t\in\R$, $h(t):=|\psi|^2_{c(t)}$ where $c$ is an arbitrary geodesic on $(\wit{M}^{2n},g)$.
We show again that $h$ is convex. 
As before $h''(t)=\mathrm{Hess}(|\psi|^2)(c'(t),c'(t))\geq 0$ for every $t\in\R$, where $\mathrm{Hess}(|\psi|^2)(X,X)=2(|p_-(X)\cdot\phi|^2+|p_+(X)\cdot\psi|^2)\geq 0$ for every $X\in T\wit{M}$ (Proposition \ref{pKirch}).
We already know that, if $V=\frac{1}{2}\grad(|\psi|^2)$ vanishes at two dif\mbox{}ferent points of $c$, then it vanishes on any geodesic segment joining the two points, therefore $|\psi|^2$ is constant on it.
This proves that $|\psi|^2$ has at most one critical value.
Since $h$ is convex this critical value is necessarily a minimum. 
The last part of the statement is a straightforward consequence of $ii)$ since $\grad(|\psi|^2)=2V$ by Proposition \ref{pKirch}.
This shows $iii)$ and concludes the proof.
\findemo


\section{Doubly warped products with imaginary K\"ahlerian Killing spinors}\label{s:dwp}

\noindent In this section, we describe the so-called doubly-warped products carrying non-zero imaginary K\"ahlerian Killing spinors.
Doubly warped products were introduced in the spinorial context by Patrick Baier in his master thesis \cite{BaierDiplA} to compute the Dirac spectrum of the complex hyperbolic space, using its representation as a doubly-warped product over an odd-dimensional sphere.\\

\noindent First we recall general formulas on warped products.

\begin{elemme}\label{l:nablawp}
Let $(\wit{M}:=M\times I,\wit{g}:=g_t\oplus\beta dt^2)$ be a warped product, where $I\subset\R$ is an open interval, $g_t$ is a smooth $1$-parameter family of Riemannian metrics on $M$ and $\beta\in C^\infty(M\times I,\R_+^\times)$.
Denote by $\wit{M}\bui{\lra}{\pi_1} M$ the first projection.
Then, for all $X,Y\in\Gamma(\pi_1^*TM)$,
\be
\wnabla_{\frac{\partial}{\partial t}}\frac{\partial}{\partial t}&=&-\frac{1}{2}\mathrm{grad}_{g_t}(\beta(t,\cdot))+\frac{1}{2\beta}\frac{\partial\beta}{\partial t}\frac{\partial}{\partial t}\\
\wnabla_{\frac{\partial}{\partial t}}X&=&\frac{\partial X}{\partial t}+\frac{1}{2}g_t^{-1}\frac{\partial g_t}{\partial t}(X,\cdot)+\frac{1}{2\beta}\frac{\partial\beta}{\partial x}(X)\frac{\partial}{\partial t}\\
\wnabla_X \frac{\partial}{\partial t}&=&\frac{1}{2}g_t^{-1}\frac{\partial g_t}{\partial t}(X,\cdot)+\frac{1}{2\beta}\frac{\partial\beta}{\partial x}(X)\frac{\partial}{\partial t}\\
\wnabla_X Y&=&\nabla_X^M Y-\frac{1}{2\beta}\frac{\partial g_t}{\partial t}(X,Y)\frac{\partial}{\partial t},
\ee
where $\frac{\partial X}{\partial t}=[\frac{\partial}{\partial t},X]$ and $\nabla^M$ (resp. $\wnabla$) is the Levi-Civita covariant derivative of $(M,g_t)$ (resp. of $(\wit{M},\wit{g})$).\\ 
\end{elemme}

\noindent{\it Proof}: straightforward consequence of the Koszul identity.
\findemo 
$ $\\

\noindent From now on we restrict ourselves to the following particular case: the manifold $M$ will be equipped with a \emph{Riemannian flow}.

\begin{edefi}\label{d:Riemflow}
\noindent\beit\item[i)]
A \emph{Riemannian flow} is a triple $(M,\wih{g},\wih{\xi})$, where $M$ is a smooth manifold and $\wih{\xi}$ is a smooth unit vector field whose flow is isometric on the orthogonal distribution, i.e., $\wih{g}(\wih{\nabla}_Z^M\wih{\xi},Z')=-\wih{g}(Z,\wih{\nabla}_{Z'}^M\wih{\xi})$ for all $Z,Z'\in\wih{\xi}^\perp$, where $\wih{\nabla}^M$ denotes the Levi-Civita covariant derivative of $(M,\wih{g})$.
\item[ii)]
A Riemannian flow $(M,\wih{g},\wih{\xi})$ is called \emph{minimal} if and only if $\wih{\nabla}_{\wih{\xi}}^M\wih{\xi}=0$, that is, if $\wih{\xi}$ is actually a Killing vector field on $M$.
\eeit
\end{edefi}
\noindent Let $(M,\wih{g},\wih{\xi})$ be a minimal Riemannian flow.
Let $\wih{h}$ denote the endomorphism-field of $\wih{\xi}^\perp$ defined by $\wih{h}(Z):=\wih{\nabla}_Z^M\wih{\xi}$ for every $Z\in\wih{\xi}^\perp$.
Let $\wih{\nabla}$ be the covariant derivative on $\wih{\xi}^\perp$ defined for all $Z\in\Gamma(\wih{\xi}^\perp)$ by $\wih{\nabla}_XZ:=\left\{\begin{array}{ll}[\wih{\xi},Z]^{\wih{\xi}^\perp}&\textrm{ if }X=\wih{\xi}\\(\wih{\nabla}_X^M Z)^{\wih{\xi}^\perp}&\textrm{ if }X\perp\wih{\xi}\end{array}\right.$.
Alternatively, $\wih{\nabla}$ can be described by the following formulas: for all $Z,Z'\in\Gamma(\wih{\xi}^\perp)$,
\[ 
\wih{\nabla}_{\wih{\xi}}^MZ=\wih{\nabla}_{\wih{\xi}} Z+\wih{h}(Z)\qquad\textrm{and}\qquad
\wih{\nabla}_Z^M Z'=\wih{\nabla}_Z Z'-\wih{g}(\wih{h}(Z),Z')\wih{\xi}.
\]
It is important to notice that, if $(M,\wih{g},\wih{\xi})$ is a (minimal) Riemannian flow and $g:=r^2(s^2\wih{g}_{\wih{\xi}}\oplus \wih{g}_{\wih{\xi}^\perp})$ for some constants $r,s>0$, then $(M,g,\xi:=\frac{1}{rs}\wih{\xi})$ is a (minimal) Riemannian flow with corresponding objects given by
\begin{equation}\label{eq:hhat}
h=\frac{s}{r}\wih{h}\qquad\textrm{and}\qquad\nabla=\wih{\nabla}.
\end{equation}

\noindent In this language, a \emph{Sasakian} manifold is a minimal Riemannian flow $(M,\wih{g},\wih{\xi})$ such that $\wih{h}$ is a transversal K\"ahler structure, that is, $\wih{h}^2=-\mathrm{Id}_{\wih{\xi}^\perp}$ and $\wih{\nabla}\wih{h}=0$.
Further on in the text we shall need for normalization purposes so-called \emph{$\mathcal{D}$-homothetic deformations} of a Sasakian structure: a $\mathcal{D}$-homothetic deformation of $(M,\wih{g},\wih{\xi})$ is $(M,\lambda^2(\lambda^2\wih{g}_{\wih{\xi}}\oplus\wih{g}_{\wih{\xi}^\perp}),\frac{1}{\lambda^2}\wih{\xi})$ for some $\lambda\in\R_+^\times$.
The identities (\ref{eq:hhat}) imply that $(M,\lambda^2(\lambda^2\wih{g}_{\wih{\xi}}\oplus\wih{g}_{\wih{\xi}^\perp}),\frac{1}{\lambda^2}\wih{\xi})$ is Sasakian as soon as $(M,\wih{g},\wih{\xi})$ is Sasakian.\\

\noindent We can now make the concept of doubly-warped product precise:

\begin{edefi}\label{d:dwp}
A \emph{doubly-warped product} is a warped product of the form
\[(\wit{M},\wit{g}):=(M\times I,\rho(t)^2(\sigma(t)^2 \wih{g}_{\wih{\xi}}\oplus \wih{g}_{\wih{\xi}^\perp})\oplus dt^2),\]
where $I$ is an open interval, $(M,\wih{g},\wih{\xi})$ is a minimal Riemannian flow, $\rho,\sigma:I\lra\R_+^\times$ are smooth functions and  $\wih{g}_{\wih{\xi}}:=\wih{g}_{|_{\R\wih{\xi}\oplus\R\wih{\xi}}}$, $\wih{g}_{\wih{\xi}^\perp}:=\wih{g}_{|_{\wih{\xi}^\perp\oplus\wih{\xi}^\perp}}$.
\end{edefi}

\noindent As for warped products, it can be easily proved that a doubly-warped product $(\wit{M},\wit{g})$ is complete as soon as $I=\R$ and $(M,\wih{g})$ is complete.\\

\noindent It is easy to check that, setting $g_t:=\rho(t)^2(\sigma(t)^2 \wih{g}_{\wih{\xi}}\oplus \wih{g}_{\wih{\xi}^\perp})$, one has $\frac{\partial g_t}{\partial t}=2\frac{\rho'}{\rho}g_t+\frac{2\sigma'}{\sigma}g_t(\pi_{\wih{\xi}^\perp},\cdot)$ and the unit vector field providing the Riemannian flow on $(M,g_t)$ is $\xi=\frac{1}{\rho\sigma}\wih{\xi}$.
In particular, the formulas in Lemma~\ref{l:nablawp} simplify:
\be  
\wnabla_{\frac{\partial}{\partial t}}\frac{\partial}{\partial t}&=&0\\
\wnabla_{\frac{\partial}{\partial t}}\xi&=&0\\
\wnabla_{\frac{\partial}{\partial t}}Z&=&\frac{\partial Z}{\partial t}+\frac{\rho'}{\rho}Z\\
\wnabla_\xi\frac{\partial}{\partial t}&=&\frac{(\rho\sigma)'}{\rho\sigma}\xi\\
\wnabla_\xi\xi&=&-\frac{(\rho\sigma)'}{\rho\sigma}\frac{\partial}{\partial t}\\
\wnabla_\xi Z&=&\nabla_\xi Z+h(Z)\\
\wnabla_Z\frac{\partial}{\partial t}&=&\frac{\rho'}{\rho}Z\\
\wnabla_Z\xi&=&h(Z)\\
\wnabla_Z Z'&=&\nabla_Z Z'-g_t(h(Z),Z')\xi-\frac{\rho'}{\rho}g_t(Z,Z')\frac{\partial}{\partial t},
\ee
where we have denoted the corresponding objects on $(M,g_t,\xi)$ without the hat ``$\widehat{\,\cdot\,}$''.\\

\noindent Next we look at a possible construction of K\"ahler structures on doubly-warped products.

\begin{elemme}\label{lstructKaehl}
Let $(\wit{M},\wit{g}):=(M\times I,\rho(t)^2(\sigma(t)^2 \wih{g}_{\wih{\xi}}\oplus \wih{g}_{\wih{\xi}^\perp})\oplus dt^2)$ be a doubly-warped product.
Assume the existence of a transversal K\"ahler structure $J$ on $(M,\wih{g},\wih{\xi})$ and define the almost complex structure $\wit{J}$ on $\wit{M}$ by
$\wit{J}(\xi):=\frac{\partial}{\partial t}$, $\wit{J}(\frac{\partial}{\partial t}):=-\xi$ and $\wit{J}(Z):=J(Z)$ for all $Z\in\{\xi,\frac{\partial}{\partial t}\}^\perp$.
Then $(\wit{M}^{2n},\wit{g},\wit{J})$ is K\"ahler if and only if $\wih{h}=-\frac{\rho'}{\sigma}J$ on $\{\xi,\frac{\partial}{\partial t}\}^\perp$ (in particular $\frac{\rho'}{\sigma}$ must be constant).
\end{elemme}

\noindent{\it Proof}: Using the identities above we write down the condition $\wnabla\wit{J}=0$.
Denote by $h$ and $\nabla$ the objects corresponding to $g_t$ on $M$.
Note first that, by definition and (\ref{eq:hhat}), one has $\nabla J=0$ on $\{\xi,\frac{\partial}{\partial t}\}^\perp$ and $\wit{J}_{|_{\{\xi,\frac{\partial}{\partial t}\}^\perp}}=J$, which does not depend on $t$.
Hence we obtain, for all $Z,Z'\in\Gamma(\wih{\xi}^\perp)$:
\be 
\wnabla_{\frac{\partial}{\partial t}}(\wit{J}(\frac{\partial}{\partial t}))-\wit{J}(\wnabla_{\frac{\partial}{\partial t}}\frac{\partial}{\partial t})&=&0\\
\wnabla_{\frac{\partial}{\partial t}}(\wit{J}(\xi))-\wit{J}(\wnabla_{\frac{\partial}{\partial t}}\xi)&=&0\\
\wnabla_{\frac{\partial}{\partial t}}(\wit{J}(Z))-\wit{J}(\wnabla_{\frac{\partial}{\partial t}}Z)&=&\frac{\partial (\wit{J}(Z))}{\partial t}-\wit{J}(\frac{\partial Z}{\partial t})=\frac{\partial (J(Z))}{\partial t}-J(\frac{\partial Z}{\partial t})=0\\
\wnabla_\xi(\wit{J}(\frac{\partial}{\partial t}))-\wit{J}(\wnabla_\xi\frac{\partial}{\partial t})&=&0\\
\wnabla_\xi(\wit{J}(\xi))-\wit{J}(\wnabla_\xi\xi)&=&0\\
\wnabla_\xi(\wit{J}(Z))-\wit{J}(\wnabla_\xi Z)&=&h\circ J(Z)-J\circ h(Z)\\
\wnabla_Z(\wit{J}(\frac{\partial}{\partial t}))-\wit{J}(\wnabla_Z\frac{\partial}{\partial t})&=&-h(Z)-\frac{\rho'}{\rho}J(Z)\\
\wnabla_Z(\wit{J}(\xi))-\wit{J}(\wnabla_Z\xi)&=&\frac{\rho'}{\rho}Z-J\circ h(Z)\\
\wnabla_Z(\wit{J}(Z'))-\wit{J}(\wnabla_Z Z')&=&-g_t(h(Z),J(Z'))\xi-\frac{\rho'}{\rho}g_t(Z,J(Z'))\frac{\partial}{\partial t}+g_t(h(Z),Z')\frac{\partial}{\partial t}-\frac{\rho'}{\rho}g_t(Z,Z')\xi.
\ee
Therefore, $\wnabla\wit{J}=0$ implies $h=-\frac{\rho'}{\rho}J$ on $\xi^\perp$ which, in turn, implies $h\circ J=J\circ h$.
Moreover, (\ref{eq:hhat}) implies that $h=\frac{\sigma}{\rho}\wih{h}$, which yields $\wih{h}=-\frac{\rho'}{\sigma}J$.
The reverse implication is obvious.
\findemo
$ $\\

\begin{erems}\label{r:witMKaehler}
{\rm\noindent\ben\item With the assumptions of Lemma~\ref{lstructKaehl}, the function $\rho'$ vanishes either identically or nowhere on the interval $I$.
In the former case the vanishing of $\wih{h}$ is equivalent to $M$ being locally the Riemannian product of an interval with a K\"ahler manifold; in the latter one, we may assume, up to changing $\sigma$ into $|\frac{\rho'}{\sigma}|\sigma$ (and $\wih{g}$ into $(\frac{\sigma}{\rho'})^2\wih{g}_{\wih{\xi}}\oplus\wih{g}_{\wih{\xi}^\perp}$), that $\wih{h}=-\varepsilon J$ and $\rho'=\varepsilon\sigma$ with $\varepsilon\in\{\pm 1\}$. 
\item Given a K\"ahler doubly warped product $(\wit{M},\wit{g},\wit{J})$ as in Lemma~\ref{lstructKaehl} and a real constant $C$, the map $(x,t)\mapsto(x,\pm t+C)$ provides a holomorphic isometry  $(\wit{M},\wit{g},\wit{J})\lra (\wit{M}',\wit{g}',\wit{J}')$, where $(\wit{M}',\wit{g}'):=(M\times (C\pm I),g_{\pm t+C}\oplus dt^2)$ and $\wit{J}'$ is the corresponding complex structure (again as in Lemma \ref{lstructKaehl}).
If furthermore $M$ is spin, then this isometry preserves the corresponding spin structures. Thus, in the case where $\rho'\neq 0$, we may assume that $\varepsilon=1$, i.e., that $\wih{h}=-J$ and $\rho'=\sigma$. 
\een
} 
\end{erems}


\noindent Now we examine the correspondence of spinors.
Let the underlying manifold $M$ of some minimal Riemannian flow $(M,g,\xi)$ be spin and, in case $M$ is the total space of a Riemannian submersion with $\S^1$-fibres over a spin manifold $N$, let $M$ carry the spin structure induced by that of $N$.
Let $\Sigma M$ denote the spinor bundle of $(M,g)$ and ``$\cdotM$'' its Clifford multiplication.
Let the doubly warped pro\-duct $\wit{M}$ carry the pro\-duct spin structure (with Clifford multiplication denoted by ``$\cdot$'').
Then the transversal covariant derivative $\nabla$ induces a covariant derivative - also denoted by $\nabla$ - on $\Sigma M$, which is related to the spinorial Levi-Civita covariant derivative $\nabla^M$ on $\Sigma M$ via (see e.g. \cite[eq. (2.4.7)]{Habibthese} or \cite[Sec. 4]{HabibEM})
\[
\nabla^M_\xi\varphi=\nabla_\xi\varphi+\frac{1}{4}\sum_{j=1}^{2n-2}e_j\cdotM h(e_j)\cdotM\varphi\qquad\textrm{and}\qquad\nabla^M_Z\varphi=\nabla_Z\varphi+\frac{1}{2}\xi\cdotM h(Z)\cdotM\varphi\]
for every $\varphi\in\Gamma(\Sigma M)$, where $\{e_j\}_{1\leq j\leq 2n-2}$ is a local orthonormal basis of $\xi^\perp\subset TM$.

\begin{elemme}\label{lspindwp}
Let a minimal Riemannian flow $(M,\wih{g},\wih{\xi})$ carry a transversal K\"ahler structure $J$ such that the doubly-warped product $(\wit{M},\wit{g},\wit{J})$ is K\"ahler, where $\wit{J}$ is the almost-complex structure induced by $J$ as in {\rm Lemma~\ref{lstructKaehl}}.
Assume furthermore $M$ to be spin.
Let $\wit{M}$ carry the induced spin structure.
Then the following identities hold for all $\varphi\in\Gamma(\Sigma\wit{M})$ and $Z\in\{\xi,\frac{\partial}{\partial t}\}^\perp$:
\be 
\wnabla_{\frac{\partial}{\partial t}}\varphi&=&\frac{\partial\varphi}{\partial t}\\
\wnabla_\xi\varphi&=&\nabla_{\xi}\varphi-\frac{\rho'}{2\rho}\wit{\Omega}\cdot\varphi-\frac{\sigma'}{2\sigma}\xi\cdot\frac{\partial}{\partial t}\cdot\varphi\\
\wnabla_Z\varphi&=&\nabla_Z\varphi-\frac{\rho'}{2\rho}(\xi\cdot J(Z)+Z\cdot\frac{\partial}{\partial t})\cdot\varphi,
\ee
where $\wit{\Omega}$ denotes the K\"ahler form of $(\wit{M},\wit{g},\wit{J})$.
\end{elemme}

\noindent{\it Proof}: Let $(e_1,\ldots,e_{2n-2},e_{2n-1}:=\xi,e_{2n}:=\frac{\partial}{\partial t})$ be a local positively-oriented orthonormal basis of $T\wit{M}$ and $(\psi_\alpha)_\alpha$ the corresponding spinorial frame.
It can be assumed that $e_j=\rho^{-1}\wih{e_j}$ with $\wih{g}(\wih{e_j},\wih{e_k})=\delta_{jk}$ and $\frac{\partial\wih{e_j}}{\partial t}=0$ (extend some $\wih{g}$-orthonormal basis independently of time).
Split $\varphi=\sum_\alpha c_\alpha\psi_\alpha$, then
\be 
\wnabla_{\frac{\partial}{\partial t}}\varphi&=&\frac{1}{4}\sum_\alpha c_\alpha\sum_{j,k=1}^{2n} \wit{g}(\wnabla_{\frac{\partial}{\partial t}} e_j,e_k)e_j\cdot e_k\cdot\psi_\alpha+\underbrace{\sum_\alpha \frac{\partial c_\alpha}{\partial t}\psi_\alpha}_{=:\frac{\partial\varphi}{\partial t}}\\
&=&\frac{\partial\varphi}{\partial t}+\frac{1}{4}\sum_\alpha c_\alpha\sum_{j,k=1}^{2n-2} \wit{g}(\wnabla_{\frac{\partial}{\partial t}} e_j,e_k)e_j\cdot e_k\cdot\psi_\alpha\\
&=&\frac{\partial\varphi}{\partial t}+\frac{1}{4}\sum_\alpha c_\alpha\sum_{j,k=1}^{2n-2}\{g_t(\frac{\partial e_j}{\partial t},e_k)+\frac{\rho'}{\rho}\delta_{jk}\}e_j\cdot e_k\cdot\psi_\alpha\\
&=&\frac{\partial\varphi}{\partial t},
\ee 
where we have used $\wnabla_{\frac{\partial}{\partial t}}\frac{\partial}{\partial t}=\wnabla_{\frac{\partial}{\partial t}}\xi=0$ and $\frac{\partial e_j}{\partial t}=-\frac{\rho'}{\rho} e_j$ by the above choice of $e_j$.
On the other hand, the Weingarten endomorphism field of $(M,g_t)$ in $\wit{M}$ is given by $A(\xi):=-\wnabla_\xi\frac{\partial}{\partial t}=-\frac{(\rho\sigma)'}{\rho\sigma}\xi$ and $A(Z):=-\wnabla_Z\frac{\partial}{\partial t}=-\frac{\rho'}{\rho}Z$ for all $Z\in\{\xi,\frac{\partial}{\partial t}\}^\perp$, so that the Gauss-Weingarten formula implies
\be 
\wnabla_\xi\varphi&=&\nabla_\xi^M\varphi+\frac{1}{2}A(\xi)\cdot\frac{\partial}{\partial t}\cdot\varphi\\
&=&\nabla_\xi\varphi+\frac{1}{4}\sum_{j=1}^{2n-2}e_j\cdotM h(e_j)\cdotM\varphi-\frac{(\rho\sigma)'}{2\rho\sigma}\xi\cdot\frac{\partial}{\partial t}\cdot\varphi\\
&=&\nabla_\xi\varphi-\frac{\rho'}{4\rho}\sum_{j=1}^{2n-2}e_j\cdot J(e_j)\cdot\varphi-\frac{(\rho\sigma)'}{2\rho\sigma}\xi\cdot\frac{\partial}{\partial t}\cdot\varphi\\
&=&\nabla_\xi\varphi-\frac{\rho'}{2\rho}\Omega\cdot\varphi-\frac{(\rho\sigma)'}{2\rho\sigma}\xi\cdot\frac{\partial}{\partial t}\cdot\varphi,
\ee
where $\Omega$ is the $2$-form associated to $J$ on $\{\xi,\frac{\partial}{\partial t}\}^\perp$, i.e., $\Omega(Z,Z')=g_t(J(Z),Z')$ for all $Z,Z'\in\{\xi,\frac{\partial}{\partial t}\}^\perp$.
Since $\wit{\Omega}=\Omega+\xi\wedge\frac{\partial}{\partial t}$, we deduce that
\be 
\wnabla_\xi\varphi&=&\nabla_\xi\varphi-\frac{\rho'}{2\rho}\wit{\Omega}\cdot\varphi+(\frac{\rho'}{2\rho}-\frac{(\rho\sigma)'}{2\rho\sigma})\xi\cdot\frac{\partial}{\partial t}\cdot\varphi\\
&=&\nabla_\xi\varphi-\frac{\rho'}{2\rho}\wit{\Omega}\cdot\varphi-\frac{\sigma'}{2\sigma}\xi\cdot\frac{\partial}{\partial t}\cdot\varphi.
\ee
For any $Z\in\{\xi,\frac{\partial}{\partial t}\}^\perp$, one has
\be 
\wnabla_Z\varphi&=&\nabla_Z^M\varphi+\frac{1}{2}A(Z)\cdot\frac{\partial}{\partial t}\cdot\varphi\\
&=&\nabla_Z\varphi+\frac{1}{2}\xi\cdotM h(Z)\cdotM\varphi-\frac{\rho'}{2\rho}Z\cdot\frac{\partial}{\partial t}\cdot\varphi\\
&=&\nabla_Z\varphi-\frac{\rho'}{2\rho}\xi\cdot J(Z)\cdot\varphi-\frac{\rho'}{2\rho}Z\cdot\frac{\partial}{\partial t}\cdot\varphi,
\ee
which shows the last identity and concludes the proof.
\findemo
$ $\\

\noindent Later on we shall need to split spinors into different components.
Recall that, on any K\"ahler spin manifold $(\wit{M}^{2n},\wit{g},\wit{J})$, the spinor bundle $\Sigma\wit{M}$ of $(\wit{M}^{2n},\wit{g})$ splits under the Clifford action of the K\"ahler form $\wit{\Omega}$ into
\[\Sigma\wit{M}=\bigoplus_{r=0}^n\Sigma_r\wit{M}, \]
where $\Sigma_r\wit{M}:=\mathrm{Ker}(\wit{\Omega}\cdot-i(2r-n)\mathrm{Id})$.
Now if $(\wit{M}^{2n},\wit{g},\wit{J})$ is a doubly-warped product as above, then any $\varphi\in\Sigma_r\wit{M}$ (with $r\in\{0,1,\ldots,n\}$) can be further split into eigenvectors for the Clifford action of $\Omega=g(J\cdot,\cdot)$.
Namely, since $[\xi\wedge\frac{\partial}{\partial t},\Omega]=0$, the automorphism $\xi\cdot\frac{\partial}{\partial t}$ of $\Sigma\wit{M}$ leaves $\Sigma_r\wit{M}$ invariant; from $(\xi\cdot\frac{\partial}{\partial t})^2=-1$ one deduces the orthogonal decomposition $\Sigma_r\wit{M}=\mathrm{Ker}(\xi\cdot\frac{\partial}{\partial t}+i\mathrm{Id})\oplus\mathrm{Ker}(\xi\cdot\frac{\partial}{\partial t}-i\mathrm{Id})$. Since both Clifford actions of $\xi$ and $\frac{\partial}{\partial t}$ are $\nabla$-parallel, so is the latter splitting.
But, for any $\varphi\in\Sigma_r\wit{M}$, one has
\be 
\varphi\in\mathrm{Ker}(\xi\cdot\frac{\partial}{\partial t}\pm i\mathrm{Id})&\iff&\Omega\cdot\varphi=i(2r-n)\varphi\pm i\varphi\\
&\iff&\Omega\cdot\varphi=i(2r-n\pm1)\varphi,
\ee
that is, $\Sigma_r\wit{M}\cap\mathrm{Ker}(\xi\cdot\frac{\partial}{\partial t}+i\mathrm{Id})=\Sigma_rM$ and $\Sigma_r\wit{M}\cap\mathrm{Ker}(\xi\cdot\frac{\partial}{\partial t}-i\mathrm{Id})=\Sigma_{r-1}M$, where by definition $\Sigma_rM:=\mathrm{Ker}(\Omega\cdot-i(2r-(n-1)\mathrm{Id}))$ for $r\in\{0,1,\ldots,n-1\}$ and $\{0\}$ otherwise.
Out of dimensional reasons one actually has
\begin{equation}\label{eqdecSr} \Sigma_r\wit{M}=\Sigma_rM\oplus\Sigma_{r-1}M\end{equation}
for every $r\in\{0,1,\ldots,n\}$.
Beware here that, if $r$ is even, then $\Sigma_r\wit{M}$ is a subspace of $\Sigma^+\wit{M}$ hence $\Sigma_r\wit{M}_{|_M}$ is canonically identified with a subspace of $\Sigma^+\wit{M}_{|_M}=\Sigma M$, whereas if $r$ is odd then it is a subspace of $\Sigma_-\wit{M}$ and is also identified as a subspace of $\Sigma M$, but this time with opposite Clifford multiplication.

\begin{elemme}\label{lspindecdwp}
Under the hypotheses of {\rm Lemma~\ref{lspindwp}}, let $\varphi\in\Gamma(\Sigma_r\wit{M})$ for some $r\in\{0,1\ldots,n\}$ and consider its decomposition $\varphi=\varphi_r+\varphi_{r-1}$ w.r.t. {\rm (\ref{eqdecSr})}. Then the identities of {\rm Lemma~\ref{lspindwp}} read:
\be 
\wnabla_{\frac{\partial}{\partial t}}\varphi_r&=&\frac{\partial\varphi_r}{\partial t}\\
\wnabla_{\frac{\partial}{\partial t}}\varphi_{r-1}&=&\frac{\partial\varphi_{r-1}}{\partial t}\\
\wnabla_\xi\varphi_r&=&\nabla_{\xi}\varphi_r+\frac{i}{2}((n-2r)\frac{\rho'}{\rho}+\frac{\sigma'}{\sigma})\varphi_r\\
\wnabla_\xi\varphi_{r-1}&=&\nabla_{\xi}\varphi_{r-1}+\frac{i}{2}((n-2r)\frac{\rho'}{\rho}-\frac{\sigma'}{\sigma})\varphi_{r-1}\\
\wnabla_Z\varphi&=&\nabla_Z\varphi_r-\frac{\rho'}{\rho}p_+(Z)\cdot\frac{\partial}{\partial t}\cdot\varphi_{r-1}+\nabla_Z\varphi_{r-1}-\frac{\rho'}{\rho}p_-(Z)\cdot\frac{\partial}{\partial t}\cdot\varphi_r
\ee
for all $Z\in\{\xi,\frac{\partial}{\partial t}\}^\perp$, where, as usual, $p_\pm(Z)=\frac{1}{2}(Z\mp iJ(Z))$.
\end{elemme}

\noindent{\it Proof}: The first two identities follow from $\wnabla_{\frac{\partial}{\partial t}}(\xi\wedge\frac{\partial}{\partial t})=0$ and $\frac{\partial J}{\partial t}=0$. For the third and fourth ones, note that $\wnabla_\xi(\xi\wedge\frac{\partial}{\partial t})=0$, so that
\be 
\wnabla_\xi\varphi_r+\wnabla_\xi\varphi_{r-1}&=&\nabla_\xi\varphi_r+\nabla_\xi\varphi_{r-1}-\frac{i\rho'}{2\rho}(2r-n)(\varphi_r+\varphi_{r-1})-\frac{i\sigma'}{2\sigma}(\varphi_{r-1}-\varphi_r)\\
&=&\nabla_\xi\varphi_r+\frac{i}{2}((n-2r)\frac{\rho'}{\rho}+\frac{\sigma'}{\sigma})\varphi_r+\nabla_\xi\varphi_{r-1}+\frac{i}{2}((n-2r)\frac{\rho'}{\rho}-\frac{\sigma'}{\sigma})\varphi_{r-1},
\ee
which is the result.
As for the last identity, one does not have $\wnabla_Z(\xi\wedge\frac{\partial}{\partial t})=0$, however 
\be 
(\xi\cdot J(Z)+Z\cdot\frac{\partial}{\partial t})\cdot\varphi&=&(-J(Z)\cdot\frac{\partial}{\partial t}\cdot\xi\cdot\frac{\partial}{\partial t}+Z\cdot\frac{\partial}{\partial t})\cdot\varphi\\
&=&-iJ(Z)\cdot\frac{\partial}{\partial t}\cdot(\varphi_{r-1}-\varphi_r)+Z\cdot\frac{\partial}{\partial t}\cdot(\varphi_r+\varphi_{r-1})\\
&=&2p_+(Z)\cdot\frac{\partial}{\partial t}\cdot\varphi_{r-1}+2p_-(Z)\cdot\frac{\partial}{\partial t}\cdot\varphi_r
\ee
for all $Z\in\{\xi,\frac{\partial}{\partial t}\}^\perp$.
This concludes the proof.
\findemo
$ $\\

\noindent We now have all we need to rewrite the imaginary K\"ahler Killing spinor equation on doubly warped products.

\begin{elemme}\label{lIKKeqdwp}
Let a spin minimal Riemannian flow $(M^{2n-1},\wih{g},\wih{\xi})$ carry a transversal K\"ahler structure $J$ such that the doubly-warped product $(\wit{M},\wit{g},\wit{J})$ is K\"ahler, where $\wit{J}$ is the almost-complex structure induced by $J$ as in {\rm Lemma~\ref{lstructKaehl}}.
Let $\wit{M}$ carry the induced spin structure and assume $n\geq 3$ to be odd.
Then a pair $(\psi,\phi)$ is an $i$-K\"ahlerian Killing spinor on $(\wit{M}^{2n},\wit{g},\wit{J})$ if and only if the following i\-den\-ti\-ties are satisfied by the components $\phi=\phi_{\frac{n+1}{2}}+\phi_{\frac{n-1}{2}}$ and $\psi=\psi_{\frac{n-1}{2}}+\psi_{\frac{n-3}{2}}$ w.r.t. {\rm (\ref{eqdecSr})}:
\begin{equation}\label{eqcharIKKSdbw} 
\left|\begin{array}{ll}\frac{\partial\phi_{\frac{n+1}{2}}}{\partial t}&=0\\
\frac{\partial\phi_{\frac{n-1}{2}}}{\partial t}&=-i\frac{\partial}{\partial t}\cdot\psi_{\frac{n-1}{2}}\\
\frac{\partial\psi_{\frac{n-1}{2}}}{\partial t}&=-i\frac{\partial}{\partial t}\cdot\phi_{\frac{n-1}{2}}\\
\frac{\partial\psi_{\frac{n-3}{2}}}{\partial t}&=0\\
\nabla_\xi\phi_{\frac{n+1}{2}}&=\frac{i}{2}(\frac{\rho'}{\rho}-\frac{\sigma'}{\sigma})\phi_{\frac{n+1}{2}}\\
\nabla_\xi\phi_{\frac{n-1}{2}}&=\frac{i}{2}(\frac{\rho'}{\rho}+\frac{\sigma'}{\sigma})\phi_{\frac{n-1}{2}}-\frac{\partial}{\partial t}\cdot\psi_{\frac{n-1}{2}}\\
\nabla_\xi\psi_{\frac{n-1}{2}}&=-\frac{i}{2}(\frac{\rho'}{\rho}+\frac{\sigma'}{\sigma})\psi_{\frac{n-1}{2}}+\frac{\partial}{\partial t}\cdot\phi_{\frac{n-1}{2}}\\
\nabla_\xi\psi_{\frac{n-3}{2}}&=-\frac{i}{2}(\frac{\rho'}{\rho}-\frac{\sigma'}{\sigma})\psi_{\frac{n-3}{2}}\\
\nabla_Z\phi_{\frac{n+1}{2}}&=p_+(Z)\cdot(\frac{\rho'}{\rho}\frac{\partial}{\partial t}\cdot\phi_{\frac{n-1}{2}}-i\psi_{\frac{n-1}{2}})\\
\nabla_Z\phi_{\frac{n-1}{2}}&=\frac{\rho'}{\rho}p_-(Z)\cdot\frac{\partial}{\partial t}\cdot\phi_{\frac{n+1}{2}}-ip_+(Z)\cdot\psi_{\frac{n-3}{2}}\\
\nabla_Z\psi_{\frac{n-1}{2}}&=\frac{\rho'}{\rho}p_+(Z)\cdot\frac{\partial}{\partial t}\cdot\psi_{\frac{n-3}{2}}-ip_-(Z)\cdot\phi_{\frac{n+1}{2}}\\
\nabla_Z\psi_{\frac{n-3}{2}}&=p_-(Z)\cdot(\frac{\rho'}{\rho}\frac{\partial}{\partial t}\cdot\psi_{\frac{n-1}{2}}-i\phi_{\frac{n-1}{2}})
\end{array}\right.\end{equation}
for every $Z\in\{\xi,\frac{\partial}{\partial t}\}^\perp$.
\end{elemme}

\noindent{\it Proof}: Since $p_+(\frac{\partial}{\partial t})\cdot\psi=\frac{1}{2}(\frac{\partial}{\partial t}+i\xi)\cdot\psi=\frac{1}{2}\frac{\partial}{\partial t}\cdot(1+i\xi\cdot\frac{\partial}{\partial t}\cdot)\psi=\frac{\partial}{\partial t}\cdot\psi_{\frac{n-1}{2}}$ and similarly $p_-(\frac{\partial}{\partial t})\cdot\phi=\frac{\partial}{\partial t}\cdot\phi_{\frac{n-1}{2}}$, the $i$-K\"ahlerian Killing spinor equation is satisfied by $(\psi,\phi)$ for $X=\frac{\partial}{\partial t}$ if and only if 
\be 
\frac{\partial\phi_{\frac{n+1}{2}}}{\partial t}+\frac{\partial\phi_{\frac{n-1}{2}}}{\partial t}&=&-ip_+(\frac{\partial}{\partial t})\cdot\psi=-i\frac{\partial}{\partial t}\cdot\psi_{\frac{n-1}{2}}\\
\frac{\partial\psi_{\frac{n-1}{2}}}{\partial t}+\frac{\partial\psi_{\frac{n-3}{2}}}{\partial t}&=&-ip_-(\frac{\partial}{\partial t})\cdot\phi=-i\frac{\partial}{\partial t}\cdot\phi_{\frac{n-1}{2}},
\ee
which gives the first four identities (use $[\Omega,\frac{\partial}{\partial t}]=0$).\\
From $p_+(\xi)\cdot\psi=-ip_+(\frac{\partial}{\partial t})\cdot\psi=-i\frac{\partial}{\partial t}\cdot\psi_{\frac{n-1}{2}}$ and $p_-(\xi)\cdot\phi=ip_-(\frac{\partial}{\partial t})\cdot\phi=i\frac{\partial}{\partial t}\cdot\phi_{\frac{n-1}{2}}$ we deduce that the $i$-K\"ahlerian Killing spinor equation is satisfied by $(\psi,\phi)$ for $X=\xi$ if and only if 
\be 
\nabla_\xi\phi_{\frac{n+1}{2}}+\frac{i}{2}(-\frac{\rho'}{\rho}+\frac{\sigma'}{\sigma})\phi_{\frac{n+1}{2}}&=&0\\
\nabla_\xi\phi_{\frac{n-1}{2}}-\frac{i}{2}(\frac{\rho'}{\rho}+\frac{\sigma'}{\sigma})\phi_{\frac{n-1}{2}}&=&-\frac{\partial}{\partial t}\cdot\psi_{\frac{n-1}{2}}\\
\nabla_\xi\psi_{\frac{n-1}{2}}+\frac{i}{2}(\frac{\rho'}{\rho}+\frac{\sigma'}{\sigma})\psi_{\frac{n-1}{2}}&=&\frac{\partial}{\partial t}\cdot\phi_{\frac{n-1}{2}}\\
\nabla_\xi\psi_{\frac{n-3}{2}}+\frac{i}{2}(\frac{\rho'}{\rho}-\frac{\sigma'}{\sigma})\psi_{\frac{n-3}{2}}&=&0,
\ee
which implies the next four equations.\\
Let $Z\in\{\xi,\frac{\partial}{\partial t}\}^\perp$, then the $i$-K\"ahlerian Killing spinor equation is satisfied by $(\psi,\phi)$ for $X=Z$ if and only if 
\be
-ip_+(Z)\cdot\psi_{\frac{n-1}{2}}&=&\nabla_Z\phi_{\frac{n+1}{2}}-\frac{\rho'}{\rho}p_+(Z)\cdot\frac{\partial}{\partial t}\cdot\phi_{\frac{n-1}{2}}\\
-ip_+(Z)\cdot\psi_{\frac{n-3}{2}}&=&\nabla_Z\phi_{\frac{n-1}{2}}-\frac{\rho'}{\rho}p_-(Z)\cdot\frac{\partial}{\partial t}\cdot\phi_{\frac{n+1}{2}}\\
-ip_-(Z)\cdot\phi_{\frac{n+1}{2}}&=&\nabla_Z\psi_{\frac{n-1}{2}}-\frac{\rho'}{\rho}p_+(Z)\cdot\frac{\partial}{\partial t}\cdot\psi_{\frac{n-3}{2}}\\
-ip_-(Z)\cdot\phi_{\frac{n-1}{2}}&=&\nabla_Z\psi_{\frac{n-3}{2}}-\frac{\rho'}{\rho}p_-(Z)\cdot\frac{\partial}{\partial t}\cdot\phi_{\frac{n-1}{2}},
\ee
which concludes the proof.
\findemo

$ $\\

\noindent Next we want to describe all doubly warped products with non-zero imaginary K\"ahlerian Killing spinors.

\begin{ethm}\label{pclassifdwp}
For $n\geq 3$ odd let $(\wit{M}^{2n},\wit{g},\wit{J})$ be a K\"ahler spin doubly warped product as in {\rm Lemma~\ref{lIKKeqdwp}}.
If there exists a non-zero $i$-K\"ahlerian Killing spinor $(\psi,\phi)$ on $(\wit{M}^{2n},\wit{g},\wit{J})$, then
\beit\item the minimal Riemannian flow $(M^{2n-1},\wih{g},\wih{\xi})$ is Sasakian,
\item up to changing $t$ into $-t$, applying a $\mathcal{D}$-homothety and translating the interval $I$ by a constant, one has either $\rho=e^t$ or $\rho=\sinh$ or $\rho=\cosh$,
\item the components $\psi_r$ and $\phi_r$ of $(\psi,\phi)$ w.r.t. {\rm (\ref{eqdecSr})} satisfy:
\eeit
\beit\item[i)] In case $\rho=e^t$: Then $\sigma=e^t$ and, setting $\wit{\psi}_{\frac{n-3}{2}}:=i\frac{\partial}{\partial t}\cdot\psi_{\frac{n-3}{2}}$ and $\varphi_{\frac{n-1}{2}}:=e^{t}(\phi_{\frac{n-1}{2}}+i\frac{\partial}{\partial t}\cdot\psi_{\frac{n-1}{2}})$, one has
\[\left|\begin{array}{ll}\frac{\partial}{\partial t}\phi_{\frac{n+1}{2}}&=0\\\frac{\partial}{\partial t}\wit{\psi}_{\frac{n-3}{2}}&=0\\\frac{\partial}{\partial t}\varphi_{\frac{n-1}{2}}&=0\\\wih{\nabla}_{\wih{\xi}}\phi_{\frac{n+1}{2}}&=0\\ \wih{\nabla}_{\wih{\xi}}\wit{\psi}_{\frac{n-3}{2}}&=0\\\wih{\nabla}\varphi_{\frac{n-1}{2}}&=0\\  \wih{\nabla}_Z\phi_{\frac{n+1}{2}}&=(-1)^{\frac{n+1}{2}}p_+(Z)\wih{\cdotM}\varphi_{\frac{n-1}{2}}\\
\wih{\nabla}_Z\wit{\psi}_{\frac{n-3}{2}}&=(-1)^{\frac{n+1}{2}}p_-(Z)\wih{\cdotM}\varphi_{\frac{n-1}{2}}.
\end{array}\right.\]
If furthermore $\varphi_{\frac{n-1}{2}}=0$, then for $\wih{\phi}_{\frac{n-1}{2}}:=e^{-t}\phi_{\frac{n-1}{2}}$ one has $\frac{\partial}{\partial t}\wih{\phi}_{\frac{n-1}{2}}=0$ and
\[\left|\begin{array}{ll}\wih{\nabla}\phi_{\frac{n+1}{2}}&=0\\\wih{\nabla}\wit{\psi}_{\frac{n-3}{2}}&=0\\\wih{\nabla}_{\wih{\xi}}\wih{\phi}_{\frac{n-1}{2}}&=0\\\wih{\nabla}_Z\wih{\phi}_{\frac{n-1}{2}}&=(-1)^{\frac{n+1}{2}}(p_-(Z)\wih{\cdotM}\phi_{\frac{n+1}{2}}+p_+(Z)\wih{\cdotM}\wit{\psi}_{\frac{n-3}{2}}).\end{array}\right.\]
In particular, the manifold $(M^{2n-1},\wih{g},\wih{\xi})$ admits a non-zero transversally parallel spinor.
Conversely, every non-zero transversally parallel spinor $\wih{\phi}_{\frac{n-1}{2}}\in\Gamma(\Sigma_{\frac{n-1}{2}}M)$ provides a non-zero $i$-K\"ahlerian Killing spinor by setting $\phi_{\frac{n+1}{2}}:=\psi_{\frac{n-3}{2}}:=0$ and $\phi_{\frac{n-1}{2}}:=e^t\wih{\phi}_{\frac{n-1}{2}}$, $\psi_{\frac{n-1}{2}}:=-e^t i\frac{\partial}{\partial t}\cdot\wih{\phi}_{\frac{n-1}{2}}$.
Moreover, for any $i$-K\"ahlerian Killing spinor $(\psi,\phi)$ on that doubly warped product $(\wit{M}^{2n},\wit{g},\wit{J})$, the component $\phi_{\frac{n-1}{2}}$ is transversally parallel on $(M,\wih{g},\wih{\xi})$ if and only if $i\frac{\partial}{\partial t}\cdot\psi=-\phi$.
\item[ii)] In case $\rho=\sinh$: One has $\sigma=\cosh$ on $I=\R_+^\times$ and there is a one-to-one correspondence between the space of $i$-K\"ahlerian Killing spinors on $(\wit{M}^{2n},\wit{g},\wit{J})$ and that of sections $(\varphi_{\frac{n+1}{2}},\varphi_{\frac{n-1}{2}},\wit{\varphi}_{\frac{n-1}{2}},\wit{\varphi}_{\frac{n-3}{2}})$ of $\Sigma_{\frac{n+1}{2}}M\oplus\Sigma_{\frac{n-1}{2}}M\oplus\Sigma_{\frac{n-1}{2}}M\oplus\Sigma_{\frac{n-3}{2}}M\lra M$ satisfying
\begin{equation}\label{eqspinKillingorig}
\left|\begin{array}{lll}\wih{\nabla}_{\wih{\xi}}\bui{\varphi}{(\sim)}_r&=&\frac{(-1)^r}{2}(n-2r)\wih{\xi}\wih{\cdotM}\bui{\varphi}{(\sim)}_r\\
\wih{\nabla}_{\wih{\xi}}\bui{\varphi}{(\sim)}_{r-1}&=&\frac{-(-1)^r}{2}(n-2r)\wih{\xi}\wih{\cdotM}\bui{\varphi}{(\sim)}_{r-1}\\
\wih{\nabla}_Z\bui{\varphi}{(\sim)}_r&=&(-1)^rp_+(Z)\wih{\cdotM}\bui{\varphi}{(\sim)}_{r-1}\\
\wih{\nabla}_Z\bui{\varphi}{(\sim)}_{r-1}&=&(-1)^rp_-(Z)\wih{\cdotM}\bui{\varphi}{(\sim)}_r
\end{array}\right.
\end{equation}
on $(M^{2n-1},\wih{g},\wih{\xi})$, for every $Z\in\wih{\xi}^\perp$ (this means that $(\varphi_{\frac{n+1}{2}},\varphi_{\frac{n-1}{2}})$ must satisfy {\rm (\ref{eqspinKillingorig})} for $r=\frac{n+1}{2}$ and $(\wit{\varphi}_{\frac{n-1}{2}},\wit{\varphi}_{\frac{n-3}{2}})$ must satisfy {\rm (\ref{eqspinKillingorig})} for $r=\frac{n-1}{2}$). 
\item[iii)] In case $\rho=\cosh$: One has $\sigma=\sinh$ on $I=\R_+^\times$ and there is a one-to-one correspondence between the space of $i$-K\"ahlerian Killing spinors on $(\wit{M}^{2n},\wit{g},\wit{J})$ and that of sections $(\varphi_{\frac{n+1}{2}},\varphi_{\frac{n-1}{2}},\wit{\varphi}_{\frac{n-1}{2}},\wit{\varphi}_{\frac{n-3}{2}})$ of $\Sigma_{\frac{n+1}{2}}M\oplus\Sigma_{\frac{n-1}{2}}M\oplus\Sigma_{\frac{n-1}{2}}M\oplus\Sigma_{\frac{n-3}{2}}M\lra M$ satisfying
\begin{equation}\label{eqspinKillinglor}
\left|\begin{array}{lll}\wih{\nabla}_{\wih{\xi}}\bui{\varphi}{(\sim)}_r&=&-\frac{(-1)^r}{2}(n-2r)\wih{\xi}\wih{\cdotM}\bui{\varphi}{(\sim)}_r\\
\wih{\nabla}_{\wih{\xi}}\bui{\varphi}{(\sim)}_{r-1}&=&\frac{(-1)^r}{2}(n-2r)\wih{\xi}\wih{\cdotM}\bui{\varphi}{(\sim)}_{r-1}\\
\wih{\nabla}_Z\bui{\varphi}{(\sim)}_r&=&(-1)^{\frac{n+1}{2}}p_+(Z)\wih{\cdotM}\bui{\varphi}{(\sim)}_{r-1}\\
\wih{\nabla}_Z\bui{\varphi}{(\sim)}_{r-1}&=&(-1)^{\frac{n-1}{2}}p_-(Z)\wih{\cdotM}\bui{\varphi}{(\sim)}_r
\end{array}\right.
\end{equation}
on $(M^{2n-1},\wih{g},\wih{\xi})$, for every $Z\in\wih{\xi}^\perp$ (this means that $(\varphi_{\frac{n+1}{2}},\varphi_{\frac{n-1}{2}})$ must satisfy {\rm (\ref{eqspinKillinglor})} for $r=\frac{n+1}{2}$ and $(\wit{\varphi}_{\frac{n-1}{2}},\wit{\varphi}_{\frac{n-3}{2}})$ must satisfy {\rm (\ref{eqspinKillinglor})} for $r=\frac{n-1}{2}$).
\eeit
\end{ethm}

\noindent{\it Proof}: We first show $\rho''=\rho$ on $I$. In order to express all equations of (\ref{eqcharIKKSdbw}) in an intrinsic way, we have to compare all objects on $(M,g_t,\xi)$ with the corresponding ones on $(M,\wih{g},\wih{\xi})$. Recall that $g_t=\rho(t)^2(\sigma(t)^2\wih{g}_{\wih{\xi}}\oplus\wih{g}_{\wih{\xi}^\perp})$ and $\xi=\frac{1}{\rho\sigma}\wih{\xi}$. 
As for (\ref{eq:hhat}), it is elementary to check the following relations:
\[\wih{\nabla}=\nabla,\;\;\;\xi\cdot=\wih{\xi}\wih{\cdot},\;\;\;\xi\cdotM=\wih{\xi}\wih{\cdotM},\;\;\;Z\cdot=\rho Z\wih{\cdot},\;\;\;Z\cdotM=\rho Z\wih{\cdotM},\]
for all $Z\in\xi^\perp$.
Applying $\frac{\partial}{\partial t}$ onto
\[\left|\begin{array}{ll}\wih{\nabla}_Z\phi_{\frac{n+1}{2}}&=p_+(Z)\wih{\cdot}(\rho'\frac{\partial}{\partial t}\cdot\phi_{\frac{n-1}{2}}-i\rho\psi_{\frac{n-1}{2}})\\ \wih{\nabla}_Z\psi_{\frac{n-3}{2}}&=p_-(Z)\wih{\cdot}(\rho'\frac{\partial}{\partial t}\cdot\psi_{\frac{n-1}{2}}-i\rho\phi_{\frac{n-1}{2}})\end{array}\right.\] and using $\frac{\partial\phi_{\frac{n+1}{2}}}{\partial t}=\frac{\partial\psi_{\frac{n-3}{2}}}{\partial t}=0$, one obtains
\be 
0&=&p_+(Z)\wih{\cdot}(\rho''\frac{\partial}{\partial t}\cdot\phi_{\frac{n-1}{2}}+\rho'\frac{\partial}{\partial t}\cdot\frac{\partial\phi_{\frac{n-1}{2}}}{\partial t}-i\rho'\psi_{\frac{n-1}{2}}-i\rho\frac{\partial\psi_{\frac{n-1}{2}}}{\partial t})\\
&=&p_+(Z)\wih{\cdot}(\rho''\frac{\partial}{\partial t}\cdot\phi_{\frac{n-1}{2}}+\rho'\frac{\partial}{\partial t}\cdot(-i\frac{\partial}{\partial t}\cdot\psi_{\frac{n-1}{2}})-i\rho'\psi_{\frac{n-1}{2}}-i\rho(-i\frac{\partial}{\partial t}\cdot\phi_{\frac{n-1}{2}}))\\
&=&(\rho''-\rho)p_+(Z)\wih{\cdot}\frac{\partial}{\partial t}\cdot\phi_{\frac{n-1}{2}}
\ee
and analogously $(\rho''-\rho)p_-(Z)\wih{\cdot}\frac{\partial}{\partial t}\cdot\psi_{\frac{n-1}{2}}=0$ for all $Z\in\wih{\xi}^\perp$.
Fix a local $\wih{g}$-orthonormal basis $(e_j)_{1\leq j\leq 2n-2}$ of $\wih{\xi}^\perp$.
Putting $Z=e_j$, Clifford-multiplying by $e_j$ and summing over $j$ gives $(\rho''-\rho)\phi_{\frac{n-1}{2}}=(\rho''-\rho)\psi_{\frac{n-1}{2}}=0$. On the other hand, both equations involving $\frac{\partial\phi_{\frac{n-1}{2}}}{\partial t}$ and $\frac{\partial\psi_{\frac{n-1}{2}}}{\partial t}$ provide the existence of smooth sections $A_{\frac{n-1}{2}}^\pm$ of $\Sigma_{\frac{n-1}{2}}M$ (independent of $t$) such that $\phi_{\frac{n-1}{2}}=e^tA_{\frac{n-1}{2}}^++e^{-t}A_{\frac{n-1}{2}}^-$ and $\psi_{\frac{n-1}{2}}=-e^t i\frac{\partial}{\partial t}\cdot A_{\frac{n-1}{2}}^++e^{-t}i\frac{\partial}{\partial t}\cdot A_{\frac{n-1}{2}}^-$.
We deduce that $(\rho''-\rho)A_{\frac{n-1}{2}}^+=(\rho''-\rho)A_{\frac{n-1}{2}}^-=0$.
If both $A_{\frac{n-1}{2}}^+$ and $A_{\frac{n-1}{2}}^-$ vanished identically on $M$, then so would $\phi_{\frac{n-1}{2}}$ and $\psi_{\frac{n-1}{2}}$ and the identities involving $\wih{\nabla}_Z\phi_{\frac{n-1}{2}}$ and $\wih{\nabla}_Z\psi_{\frac{n-1}{2}}$ would provide (after contracting with the Clifford multiplication just as above)  $\phi_{\frac{n+1}{2}}=\psi_{\frac{n-3}{2}}=0$, so that $(\psi,\phi)=0$, which is a contradiction.
Therefore $\rho''-\rho=0$ on $I$.\\
It follows in particular that $\rho'=0$ on $I$ cannot hold, so we may assume that $\wih{h}=- J$ (hence $(M^{2n-1},\wih{g},\wih{\xi})$ is Sasakian) and $\rho'=\sigma$ (see Remarks \ref{r:witMKaehler}).
Furthermore, in the case where the constant $(\rho')^2-\rho^2$ does not vanish, up to replacing $\rho$ by $\frac{\rho}{\sqrt{|(\rho'^2)-\rho^2|}}$ (which is equivalent to performing a $\mathcal{D}$-homothetic deformation of the Sasakian structure), we may assume that $(\rho'^2)-\rho^2=1$ or $-1$ on $I$.
Next we rewrite the equations from Lemma~\ref{lIKKeqdwp} considering the new sections $\varphi_{\frac{n+1}{2}},\varphi_{\frac{n-1}{2}},\wit{\varphi}_{\frac{n-1}{2}},\wit{\varphi}_{\frac{n-3}{2}}$ defined by
\[\left|\begin{array}{ll}\varphi_{\frac{n+1}{2}}&:=\phi_{\frac{n+1}{2}}\\
\varphi_{\frac{n-1}{2}}&:=\rho'\phi_{\frac{n-1}{2}}+i\rho\frac{\partial}{\partial t}\cdot\psi_{\frac{n-1}{2}}\\
\wit{\varphi}_{\frac{n-1}{2}}&:=i\rho\frac{\partial}{\partial t}\cdot\phi_{\frac{n-1}{2}}+\rho'\psi_{\frac{n-1}{2}}\\
\wit{\varphi}_{\frac{n-3}{2}}&:=\psi_{\frac{n-3}{2}}.\end{array}\right.\]
Note that the linear transformation $(\phi_{\frac{n+1}{2}},\phi_{\frac{n-1}{2}},\psi_{\frac{n-1}{2}},\psi_{\frac{n-3}{2}})\mapsto(\varphi_{\frac{n+1}{2}},\varphi_{\frac{n-1}{2}},\wit{\varphi}_{\frac{n-1}{2}},\wit{\varphi}_{\frac{n-3}{2}})$ is invertible if and only if $(\rho')^2-\rho^2\neq0$. From (\ref{eqcharIKKSdbw}) we have, for all $Z\in\wih{\xi}^\perp$:
\be 
\frac{\partial}{\partial t}\varphi_{\frac{n+1}{2}}&=&0\\
\frac{\partial}{\partial t}\varphi_{\frac{n-1}{2}}&=&0\\
\frac{\partial}{\partial t}\wit{\varphi}_{\frac{n-1}{2}}&=&0\\
\frac{\partial}{\partial t}\wit{\varphi}_{\frac{n-3}{2}}&=&0\\
\wih{\nabla}_{\wih{\xi}}\varphi_{\frac{n+1}{2}}&=&\frac{(-1)^{\frac{n+1}{2}}}{2}(n-2(\frac{n+1}{2}))((\rho')^2-\rho^2)\wih{\xi}\wih{\cdotM}\varphi_{\frac{n+1}{2}}\\
\wih{\nabla}_{\wih{\xi}}\varphi_{\frac{n-1}{2}}&=&-\frac{(-1)^{\frac{n+1}{2}}}{2}(n-2(\frac{n+1}{2}))((\rho')^2-\rho^2)\wih{\xi}\wih{\cdotM}\varphi_{\frac{n-1}{2}}\\
\wih{\nabla}_{\wih{\xi}}\wit{\varphi}_{\frac{n-1}{2}}&=&\frac{(-1)^{\frac{n-1}{2}}}{2}(n-2(\frac{n-1}{2}))((\rho')^2-\rho^2)\wih{\xi}\wih{\cdotM}\wit{\varphi}_{\frac{n-1}{2}}\\
\wih{\nabla}_{\wih{\xi}}\wit{\varphi}_{\frac{n-3}{2}}&=&-\frac{(-1)^{\frac{n-1}{2}}}{2}(n-2(\frac{n-1}{2}))((\rho')^2-\rho^2)\wih{\xi}\wih{\cdotM}\wit{\varphi}_{\frac{n-3}{2}}\\
\wih{\nabla}_Z\varphi_{\frac{n+1}{2}}&=&(-1)^{\frac{n+1}{2}}p_+(Z)\wih{\cdotM}\varphi_{\frac{n-1}{2}}\\
\wih{\nabla}_Z\varphi_{\frac{n-1}{2}}&=&(-1)^{\frac{n+1}{2}}((\rho')^2-\rho^2)p_-(Z)\wih{\cdotM}\varphi_{\frac{n+1}{2}}\\
\wih{\nabla}_Z\wit{\varphi}_{\frac{n-1}{2}}&=&(-1)^{\frac{n-1}{2}}((\rho')^2-\rho^2)p_+(Z)\wih{\cdotM}\wit{\varphi}_{\frac{n-3}{2}}\\
\wih{\nabla}_Z\wit{\varphi}_{\frac{n-3}{2}}&=&(-1)^{\frac{n-1}{2}}p_-(Z)\wih{\cdotM}\wit{\varphi}_{\frac{n-1}{2}}.\\
\ee
If $(\rho')^2-\rho^2\neq 0$ on $I$, then the required equations directly follow from the above ones. Moreover, since in that case the correspondence $(\phi_{\frac{n+1}{2}},\phi_{\frac{n-1}{2}},\psi_{\frac{n-1}{2}},\psi_{\frac{n-3}{2}})\mapsto(\varphi_{\frac{n+1}{2}},\varphi_{\frac{n-1}{2}},\wit{\varphi}_{\frac{n-1}{2}},\wit{\varphi}_{\frac{n-3}{2}})$ is bijective, the ``If'' in the assumptions is actually an ``if and only if''.
If now $(\rho')^2-\rho^2=0$, then $\rho'=\pm\rho$ on $I$; since we have assumed $\rho'>0$ (up to changing $t$ into $-t$), we only have to consider $\rho'=\rho$, hence $\rho=Ce^t$ for some positive constant $C$.
Since translating $t$ provides a holomorphic isometry (again see Remarks~\ref{r:witMKaehler}), one may assume that $C=1$, i.e., $\rho=e^t$.
In that case, one has $\wih{\nabla}\varphi_{\frac{n-1}{2}}=0$ on $M$, hence $\varphi_{\frac{n-1}{2}}$ vanishes either identically or nowhere on $M$ (and on $\wit{M}$ since it is constant in $t$).
If $\varphi_{\frac{n-1}{2}}\neq 0$, then all right members in the equations listed just above vanish except 
\[ \left|\begin{array}{ll}\wih{\nabla}_Z\phi_{\frac{n+1}{2}}&=(-1)^{\frac{n+1}{2}}p_+(Z)\wih{\cdotM}\varphi_{\frac{n-1}{2}}\\ 
\wih{\nabla}_Z\wit{\varphi}_{\frac{n-3}{2}}&=(-1)^{\frac{n-1}{2}}p_-(Z)\wih{\cdotM}\wit{\varphi}_{\frac{n-1}{2}},\end{array}\right.\] 
which together with $\wit{\varphi}_{\frac{n-1}{2}}=i\frac{\partial}{\partial t}\cdot \varphi_{\frac{n-1}{2}}$ gives the result.
If $\varphi_{\frac{n-1}{2}}=0$ on $M$, then coming back to the equations from Lemma~\ref{lIKKeqdwp}, one has $\wih{\nabla}\phi_{\frac{n+1}{2}}=\wih{\nabla}\psi_{\frac{n-3}{2}}=0$ and $\wih{\phi}_{\frac{n-1}{2}}$ satisfies the required equations.
\findemo  
$ $\\

\begin{erem}\label{rexspKKI}
{\rm In Theorem \ref{pclassifdwp}$.i)$ not every $i$-K\"ahlerian Killing spinor on $\wit{M}$ must come from a transversally parallel spinor on $M$.
For instance, consider the complex hyperbolic space $\C\mathrm{H}^n$ (for $n$ odd) endowed with its Fubini-Study metric of constant holomorphic sectional curvature $-4$ and its canonical spin structure.
Then $\C\mathrm{H}^n$ (possibly with a suitable submanifold removed) can be viewed as a doubly warped product in several ways.
For example, $\C\mathrm{H}^n$ is a doubly-warped product over the Heisenberg group $M$, which admits a $\left(\begin{array}{c}n-1\\\frac{n-1}{2}\end{array}\right)$-dimensional space of transversally parallel spinors lying pointwise in $\Sigma_{\frac{n-1}{2}}M$ (see below).
However, $\C\mathrm{H}^n$ carries a $2\left(\begin{array}{c}n\\\frac{n+1}{2}\end{array}\right)$-dimensional space of $i$-K\"ahlerian Killing spinors \cite[Sec. 3]{KirchbergAGAG93}.
Therefore there exists at least one non-zero K\"ahlerian Killing spinor on $\C\mathrm{H}^n$ which does not come from any transversally parallel spinor on $M$.}
\end{erem}

\noindent As an example for Theorem~\ref{pclassifdwp}$.i)$, any Heisenberg manifold of dimension $4k+1$ (with $k\geq 1$) has a spin structure for which the corresponding spinor bundle is trivialized by transversally parallel spinors. 
This follows from three facts: every Heisenberg manifold is an $\S^1$-bundle with totally geodesic fibres over a flat torus; every $\S^1$-bundle over a manifold carrying parallel spinors carries transversally parallel spinors for the induced spin structure, see e.g. \cite[Prop. 3.6]{GinHabibspKTgeom}; the whole spinor bundle of any flat torus endowed with its so-called trivial spin structure is trivialized by parallel spinors.
Note that, as a consequence of Lemma~\ref{lholseccurv} below, the doubly warped product arising from a $(2n-1)$-dimensional Heisenberg manifold $M$ choosing $\rho=\sigma=e^t$ has constant holomorphic sectional curvature $-4$, therefore it is holomorphically isometric to $\C\mathrm{H}^n$ as soon as it is simply-connected and complete.\\
Examples for Theorem \ref{pclassifdwp}$.i)$ with non-constant holomorphic sectional curvature can be constructed out of the following lemma:

\begin{elemme}\label{ltransvpar}
For each integer $n\equiv 1\;(4)$, let $(N^{2n-2},g_N,J)$ be any simply-connected closed Hodge hyperk\"ahler manifold.
Then there exists an $\S^1$-bundle $M$ over $N$ carrying an $\S^1$-invariant metric $\wih{g}$ for which $(M^{2n-1},\wih{g},\wih{\xi})$ is Sasakian and for which there exists a parallel spinor lying pointwise in $\Sigma_{\frac{n-1}{2}}M$. 
\end{elemme}

\noindent{\it Proof}: Recall first that every hyperk\"ahler manifold is spin (this follows from the structure group $\mathrm{Sp}(\frac{n-1}{2})$ being simply-connected).
McK. Wang's classification \cite{Wang89} of manifolds with parallel spinors provides the existence of exactly $\frac{n-1}{2}+1$ linearly independent parallel spinors on $N$, one of which lies pointwise in $\Sigma_{\frac{n-1}{2}}N$ if and only if $\frac{n-1}{2}$ is even \cite[(ii) p.61]{Wang89}.
Now, for any Hodge K\"ahler manifold $(N,g,J)$ (``Hodge'' meaning that its K\"ahler class is proportional to an integral class), there exists an $\S^1$-bundle $M\bui{\lra}{\pi}N$ carrying an $\S^1$-invariant metric $\wih{g}$ for which $(M^{2n-1},\wih{g},\wih{\xi})$ is Sasakian with $\wih{h}=-J$, see \cite[Prop. 2]{MoroiHodge} (as usual $\wih{\xi}$ denotes the fundamental vector field of the $\S^1$-action).
By \cite[Prop. 3.6]{GinHabibspKTgeom}, the lift of the non-zero parallel spinor in $\Sigma_{\frac{n-1}{2}}N$ to $M$ gives a non-zero transversal parallel spinor on $(M^{2n-1},\wih{g},\wih{\xi})$ provided the spin structure on $M$ is induced by the one on $\pi^*(TN)$ and the trivial covering of $\S^1$; because of $\wih{h}=-J$, this spinor lies pointwise in $\Sigma_{\frac{n-1}{2}}M$.
\findemo
 $ $\\

\noindent Kodaira's embedding theorem states that a closed K\"ahler manifold is Hodge if and only if it is projective, i.e., if and only if it can be holomorphically embedded in some complex projective space.
Therefore projective hyperk\"ahler manifolds of complex dimension $4k$ (with $k\geq 1$) provide examples for $N$ in Lemma~\ref{ltransvpar}.
For instance, simply connected hyperk\"ahler manifolds can be constructed as
the Hilbert scheme of a ${\mathrm K3}$-surface (cf. \cite{Beauville83}).
Indeed, let $X$ be a ${\mathrm K3}$-surface, then the Hilbert scheme
$\mathrm{Hilb}^{2k}(X)$, which is the blow-up along the diagonal of the $2k$-th
symmetric product of $X$,  is a compact, simply-connected hyperk\"ahler
manifold of complex dimension $4k$.
If $X$ is projective, e.g. a quartic,
then $\mathrm{Hilb}^{2k}(X)$ is projective too and thus has an
integer K\"ahler class.\\

\noindent In order to decide whether the doubly warped product we construct is the complex hyperbolic space or not, the transversal holomorphic curvature of $(M,\wih{g},\wih{\xi})$ and the holomorphic sectional curvature of $(\wit{M}^{2n},\wit{g},\wit{J})$ have to be compared:

\begin{elemme}\label{lholseccurv}
Let $(\wit{M}^{2n},\wit{g},\wit{J})$ be a K\"ahler doubly warped product as in {\rm Lemma~\ref{lstructKaehl}} with $\rho''=\rho$, $\sigma=\rho'$ and $\wih{h}=-J$.
Then the holomorphic sectional curvature $\wit{K}_{\rm hol}(Z)$ of $(\wit{M},\wit{g},J)$ and the transversal holomorphic sectional curvature $\wih{K}_{\rm hol}(Z)$ of $(M,\wih{g},\wih{\xi})$ are related by
\[\wit{K}_{\rm hol}(Z)=\frac{1}{\rho^2}\Big(\wih{K}_{\rm hol}(Z)-4(\rho')^2\Big),\]
for all $Z\in\{\wih{\xi},\frac{\partial}{\partial t}\}^\perp\setminus\{0\}$.
In particular, the doubly warped product $(\wit{M}^{2n},\wit{g},\wit{J})$ has constant holomorphic sectional curvature $-4$ if and only if the transversal holomorphic sectional curvature of $(M,\wih{g},\wih{\xi})$ is constant equal to $4((\rho')^2-\rho^2)$.
\end{elemme}

\noindent{\it Proof}: Recall that $\wit{K}_{\rm hol}(Z)$ and $\wih{K}_{\rm hol}(Z)$ are defined by
\[\wit{K}_{\rm hol}(Z):=\frac{\wit{g}(\wit{R}(Z,JZ)Z,JZ)}{\wit{g}(Z,Z)^2}\qquad\textrm{and}\qquad\wih{K}_{\rm hol}(Z):=\frac{\wih{g}(\wih{R}(Z,JZ)Z,JZ)}{\wih{g}(Z,Z)^2},\]
where $\wit{R}_{X,Y}:=\wnabla_{[X,Y]}-[\wnabla_X,\wnabla_Y]$ and $\wih{R}_{Z,Z'}:=\wih{\nabla}_{[Z,Z']}-[\wih{\nabla}_Z,\wih{\nabla}_{Z'}]$ are the curvature tensors associated to $\wnabla$ and $\wih{\nabla}$ on $T\wit{M}$ and $\wih{\xi}^\perp$ respectively. The following identities can be deduced from the formulas in Lemma~\ref{l:nablawp}, taking into account $\rho'=\sigma$ and $\rho''=\rho$:
\be 
\wit{g}(\wit{R}_{\xi,\frac{\partial}{\partial t}}\xi,\frac{\partial}{\partial t})&=&-\frac{(\rho\sigma)''}{\rho\sigma}=-4\\
\wit{g}(\wit{R}(Z,JZ)Z,JZ)&=&\wit{g}(\wih{R}(Z,JZ)Z,JZ)-4(\frac{\rho'}{\rho})^2\wit{g}(Z,Z)^2,
\ee
for every $Z\in\{\wih{\xi},\frac{\partial}{\partial t}\}^\perp\setminus\{0\}$.
Using $\wit{g}(Z,\cdot)=\rho^2\wih{g}(Z,\cdot)$, we obtain
\be 
\wit{K}_{\rm hol}(Z)&=&\frac{\wit{g}(\wih{R}(Z,JZ)Z,JZ)}{\wit{g}(Z,Z)^2}-4(\frac{\rho'}{\rho})^2\\
&=&\frac{1}{\rho^2}\frac{\wih{g}(\wih{R}(Z,JZ)Z,JZ)}{\wih{g}(Z,Z)^2}-4(\frac{\rho'}{\rho})^2,
\ee
which gives the first statement.
Since by the computation above $\wit{K}_{\rm hol}(\xi)=-4$ (independently of $\wih{g}$), the second follows from the first (note that $(\rho')^2-\rho^2$ is constant by the assumption $\rho''=\rho$).
\findemo
$ $\\

\noindent As a consequence of Theorem~\ref{pclassifdwp}$.i)$,  Lemma~\ref{ltransvpar} and Lemma~\ref{lholseccurv}, we obtain:

\begin{ecor}\label{c:exnoncstholsect}
For an integer $n\equiv 1\;(4)$, let $(N^{2n-2},g_N,J)$ be any simply-connected closed Hodge hyperk\"ahler manifold.
Let $(M^{2n-1},\wih{g},\wih{\xi})$ be constructed from $N$ as in {\rm Lemma \ref{ltransvpar}} and $(\wit{M}^{2n},\wit{g},\wit{J})$ be the K\"ahler spin doubly warped product constructed from $M$ as in {\rm Lemma~\ref{lspindwp}} with $\rho=\sigma=e^t$.
Then $(\wit{M}^{2n},\wit{g},\wit{J})$ carries a non-zero $i$-K\"ahlerian Killing spinor but has non-constant holomorphic sectional curvature.
\end{ecor}

\noindent{\it Proof}: The existence of a non-zero $i$-K\"ahlerian Killing spinor follows from Theorem~\ref{pclassifdwp}$.i)$ and Lemma~\ref{ltransvpar}.
In case $\rho=\sigma=e^t$, Lemma~\ref{lholseccurv} implies that the holomorphic sectional curvature of the doubly warped product $(\wit{M}^{2n},\wit{g},J)$ is $-4$ if and only if the transversal holomorphic sectional curvature of $(M,\wih{g},\wih{\xi})$ vanishes, that is, if and only if its transversal curvature vanishes (see e.g. \cite[Prop. 7.1 p.166]{KoNo}).
Now for any $\S^1$-bundle as in Lemma~\ref{ltransvpar}, the transversal (holomorphic) sectional curvature of $M$ and the (holomorphic) sectional curvature of $N$ coincide.
Since simply-connected closed hyperk\"ahler manifolds cannot be flat, the K\"ahler manifold $(\wit{M}^{2n},\wit{g},\wit{J})$ cannot have constant holomorphic sectional curvature.
\findemo
$ $\\

\noindent Corollary~\ref{c:exnoncstholsect} provides the first family of examples of K\"ahler spin manifolds of non-constant holomorphic sectional curvature carrying non-zero imaginary K\"ahlerian Killing spinors.\\

\noindent The two other subcases $(\rho')^2-\rho^2=1$ and $(\rho')^2-\rho^2=-1$ are geometrically more simple to describe.
We do it in separate lemmas.

\begin{elemme}\label{lequivspinKilling}
Let $(M^{2n-1},g,\xi)$ be a Sasakian spin manifold with $h=-J$ and fix $r\in\{0,1,\ldots,n\}$. Then a section $(\psi_r,\psi_{r-1})$ of $\Sigma_rM\oplus\Sigma_{r-1}M$ satisfies {\rm (\ref{eqspinKillingorig})} if and only if $\psi:=\psi_r+\psi_{r-1}$ is a $\frac{(-1)^r}{2}$-Killing spinor on $(M,g)$.
\end{elemme}

\noindent{\it Proof}: Let $\Omega$ be the $2$-form associated to $J$ on $\xi^\perp$, i.e., $\Omega(Z,Z')=g(J(Z),Z')$ for all $Z,Z'\perp\xi$. Using $\Omega\cdotM\psi_r=(-1)^{r+1}(2r-n+1)\xi\cdotM\psi_r$ (for all $r$) we have on the one hand
\be
\nabla_\xi\psi&=&\nabla_\xi^M\psi+\frac{1}{2}\Omega\cdotM\psi\\
&=&\nabla_\xi^M\psi-\frac{(-1)^r}{2}\xi\cdotM\psi+\frac{(-1)^r}{2}\xi\cdotM\psi+\frac{1}{2}\Omega\cdotM\psi\\
&=&\nabla_\xi^M\psi-\frac{(-1)^r}{2}\xi\cdotM\psi+\frac{(-1)^r}{2}\xi\cdotM\psi-\frac{(-1)^r}{2}(2r-n+1)\xi\cdotM\psi_r\\
& &+\frac{(-1)^r}{2}(2(r-1)-n+1)\xi\cdotM\psi_{r-1}\\
&=&\nabla_\xi^M\psi-\frac{(-1)^r}{2}\xi\cdotM\psi+\frac{(-1)^r}{2}(n-2r)\xi\cdotM\psi_r+\frac{(-1)^r}{2}(2(r-1)-n+2)\xi\cdotM\psi_{r-1},
\ee
which implies
\begin{equation}\label{eqxi}
\left|\begin{array}{ll}\nabla_\xi\psi_r&=(\nabla_\xi^M\psi-\frac{(-1)^r}{2}\xi\cdotM\psi)_r+\frac{(-1)^r}{2}(n-2r)\xi\cdotM\psi_r\\
 & \\
\nabla_\xi\psi_{r-1}&=(\nabla_\xi^M\psi-\frac{(-1)^r}{2}\xi\cdotM\psi)_{r-1}-\frac{(-1)^r}{2}(n-2r)\xi\cdotM\psi_{r-1}.
\end{array}\right.
\end{equation}
On the other hand, for every $Z\in\xi^\perp$ one has, 
\be
\nabla_Z\psi&=&\nabla_Z^M\psi-\frac{1}{2}\xi\cdotM h(Z)\cdotM\psi\\
&=&\nabla_Z^M\psi-\frac{(-1)^r}{2}Z\cdotM\psi+\frac{(-1)^r}{2}Z\cdotM\psi-\frac{1}{2}J(Z)\cdotM\xi\cdotM\psi\\
&=&\nabla_Z^M\psi-\frac{(-1)^r}{2}Z\cdotM\psi+\frac{(-1)^r}{2}Z\cdotM\psi-\frac{1}{2}J(Z)\cdotM\{(-1)^{r+1}i\psi_r+(-1)^r i\psi_{r-1}\}\\
&=&\nabla_Z^M\psi-\frac{(-1)^r}{2}Z\cdotM\psi+\frac{(-1)^r}{2}(Z+iJ(Z))\cdot\psi_r+\frac{(-1)^r}{2}(Z-iJ(Z))\cdot\psi_{r-1},
\ee
which implies
\begin{equation}\label{eqX}
\left|\begin{array}{ll}\nabla_Z\psi_r&=(\nabla_Z^M\psi-\frac{(-1)^r}{2}Z\cdotM\psi)_r+(-1)^rp_+(Z)\cdot\psi_{r-1}\\
\nabla_Z\psi_{r-1}&=(\nabla_Z^M\psi-\frac{(-1)^r}{2}Z\cdotM\psi)_{r-1}+(-1)^rp_-(Z)\cdot\psi_r.\end{array}\right.
\end{equation}
Therefore the pair $(\psi_r,\psi_{r-1})$ satisfies (\ref{eqspinKillingorig}) if and only if $\psi:=\psi_r+\psi_{r-1}$ satisfies $\nabla_X^M\psi=\frac{(-1)^r}{2}X\cdotM\psi$ for all $X\in TM$, that is, if and only if $\psi$ is a $\frac{(-1)^r}{2}$-Killing spinor on $(M,g)$.
\findemo
$ $\\

\noindent The case $(\rho')^2-\rho^2=-1$ is analogous to the case $(\rho')^2-\rho^2=1$ up to a Lorentzian detour. We call {\rm (\ref{eqspinKillinglorr})} the following system of equations:
\begin{equation}\label{eqspinKillinglorr}
\left|\begin{array}{ll} 
\nabla_\xi\psi_r&=-\frac{(-1)^r}{2}(n-2r)\xi\cdotM\psi_r\\
\nabla_\xi\psi_{r-1}&=\frac{(-1)^r}{2}(n-2r)\xi\cdotM\psi_{r-1}\\
\nabla_Z\psi_r&=(-1)^r\epsilon p_+(Z)\cdotM\psi_{r-1}\\
\nabla_Z\psi_{r-1}&=-(-1)^r\epsilon p_-(Z)\cdotM\psi_{r}
\end{array}\right.
\end{equation}
for all $Z,Z'\in\xi^\perp$, where $\epsilon\in\{\pm1\}$.

\begin{elemme}\label{lequivspinKillinglor}
Let $(M^{2n-1},g,\xi)$ be a Sasakian spin manifold with $h=-J$ and fix $r\in\{0,1,\ldots,n\}$ as well as $\epsilon\in\{\pm1\}$. Then a section $(\psi_r,\psi_{r-1})$ of $\Sigma_rM\oplus\Sigma_{r-1}M$ satisfies {\rm (\ref{eqspinKillinglorr})} if and only if $\psi:=\psi_r+i\epsilon\psi_{r-1}$ is a $\frac{(-1)^{r+1} i}{2}$-Killing spinor on the Lorentzian manifold $(M,-g_{\xi}\oplus g_{\xi^\perp})$.
\end{elemme}

\noindent{\it Proof}: First, there exists the analog of Riemannian flow in the Lorentzian context. A Lorentzian flow is given by a triple $(M,\wih{g},\wih{\xi})$, where $(M,\wih{g})$ is a Lorentzian manifold and $\wih{\xi}$ a smooth tangent vector field on $M$ with $\wih{g}(\wih{\xi},\wih{\xi})=-1$ and $\wih{g}(\wih{\nabla}_Z^M\wih{\xi},Z')=-\wih{g}(\wih{\nabla}_{Z'}^M\wih{\xi},Z)$ for all $Z,Z'\in\wih{\xi}^\perp$. Note that $(M,\wih{g})$ is necessarily time-oriented because of the existence of $\wih{\xi}$. Setting $\wih{\nabla}_XZ:=\left|\begin{array}{ll}[\wih{\xi},Z]^{\wih{\xi}^\perp}&\textrm{ if }X=\wih{\xi}\\(\wih{\nabla}_X^M Z)^{\wih{\xi}^\perp}&\textrm{ if }X\perp\wih{\xi}\end{array}\right.$ for all $Z\in\Gamma(\wih{\xi}^\perp)$ and $\wih{h}:=\wih{\nabla}^M\wih{\xi}$, one obtains a metric connection $\wih{\nabla}$ and a skew-symmetric endomorphism-field $\wih{h}$ on $\wih{\xi}^\perp$ such that
\[ \left|\begin{array}{ll}
\wih{\nabla}_{\wih{\xi}}^MZ&=\wih{\nabla}_{\wih{\xi}}Z+\wih{h}(Z)+\wih{g}(\wih{\nabla}_{\wih{\xi}}^M\wih{\xi},Z)\wih{\xi}\\
\wih{\nabla}_{Z}^MZ'&=\wih{\nabla}_{Z}Z'+\wih{g}(\wih{h}(Z),Z')\wih{\xi}
\end{array}\right.\]
for all $Z,Z'\in\Gamma(\wih{\xi}^\perp)$. Moreover, in case $M$ is spin, the corresponding Gauss-type formula for spinors reads
\[\left|\begin{array}{ll}
\wih{\nabla}_{\wih{\xi}}\varphi&=\wih{\nabla}_{\wih{\xi}}^M\varphi-\frac{1}{2}\wih{\Omega}\wih{\cdotM}\varphi+\frac{1}{2}\wih{\xi}\wih{\cdotM}\wih{\nabla}_{\wih{\xi}}^M\wih{\xi}\wih{\cdotM}\varphi\\ 
\wih{\nabla}_Z\varphi&=\wih{\nabla}_Z^M\varphi+\frac{1}{2}\wih{\xi}\wih{\cdotM}\wih{h}(Z)\wih{\cdotM}\varphi 
\end{array}\right. \]
for all $\varphi\in\Gamma(\Sigma M)$ and $Z\in\wih{\xi}^\perp$, where $\wih{\Omega}(Z,Z'):=\wih{g}(\wih{h}(Z),Z')$. In case $(M,\wih{g},\wih{\xi})$ is Lorentzian Sasakian, i.e., if furthermore $\wih{\nabla}_{\wih{\xi}}^M\wih{\xi}=0$, $\wih{h}^2=-\mathrm{Id}$ and $\wih{\nabla}\wih{h}=0$, then we still have the $\wih{\nabla}$-parallel decomposition $\Sigma M=\oplus_{r=0}^{n-1}\Sigma_rM$ with $\Sigma_rM:=\mathrm{Ker}(\wih{\Omega}\wih{\cdotM}-i(2r-(n-1)\mathrm{Id}))$. This time one has $\wih{\xi}\wih{\cdotM}\varphi_r=(-1)^{r+1}\varphi_r$ for all $\varphi_r\in\Sigma_rM$.\\
Assume now $(M,\wih{g},\wih{\xi})$ to be Lorentzian Sasakian and pick a section $\psi=\psi_r+\psi_{r-1}$ of $\Sigma_r M\oplus\Sigma_{r-1} M$, then the formulas above imply
\be 
\wih{\nabla}_{\wih{\xi}}\psi&=&\wih{\nabla}_{\wih{\xi}}^M\psi-\frac{1}{2}\wih{\Omega}\wih{\cdotM}\psi\\
&=&\wih{\nabla}_{\wih{\xi}}^M\psi-\frac{(-1)^{r+1} i}{2}\wih{\xi}\wih{\cdotM}\psi+\frac{(-1)^{r+1} i}{2}\wih{\xi}\wih{\cdotM}\psi-\frac{i}{2}\big((2r-(n-1))\psi_r+(2(r-1)-(n-1))\psi_{r-1}\big)\\
&=&\wih{\nabla}_{\wih{\xi}}^M\psi-\frac{(-1)^{r+1} i}{2}\wih{\xi}\wih{\cdotM}\psi+\frac{(-1)^{r+1} i}{2}\wih{\xi}\wih{\cdotM}\psi\\
& &+\frac{(-1)^{r}i}{2}(2r-(n-1))\wih{\xi}\wih{\cdotM}\psi_r-\frac{(-1)^r i}{2}(2(r-1)-(n-1))\wih{\xi}\wih{\cdotM}\psi_{r-1}\\
&=&\wih{\nabla}_{\wih{\xi}}^M\psi-\frac{(-1)^{r+1} i}{2}\wih{\xi}\wih{\cdotM}\psi+\frac{(-1)^{r+1} i}{2}(n-2r)\wih{\xi}\wih{\cdotM}\psi_r-\frac{(-1)^{r+1} i}{2}(n-2r)\wih{\xi}\wih{\cdotM}\psi_{r-1},
\ee
that is,
\[\left|\begin{array}{ll}\wih{\nabla}_{\wih{\xi}}\psi_r&=\big(\wih{\nabla}_{\wih{\xi}}^M\psi-\frac{(-1)^{r+1} i}{2}\wih{\xi}\wih{\cdotM}\psi\big)_r+\frac{(-1)^{r+1} i}{2}(n-2r)\wih{\xi}\wih{\cdotM}\psi_r\\ 
\wih{\nabla}_{\wih{\xi}}\psi_{r-1}&=\big(\wih{\nabla}_{\wih{\xi}}^M\psi-\frac{(-1)^{r+1} i}{2}\wih{\xi}\wih{\cdotM}\psi\big)_{r-1}-\frac{(-1)^{r+1} i}{2}(n-2r)\wih{\xi}\wih{\cdotM}\psi_{r-1}.\end{array}\right.\]
This is still valid for $r=0$ or $r=n$ (setting $\psi_{-1}:=\psi_n:=0$). Similarly, for all $Z\in\wih{\xi}^\perp$,
\be 
\wih{\nabla}_Z\psi&=&\wih{\nabla}_Z^M\psi+\frac{1}{2}\wih{\xi}\wih{\cdotM}\wih{h}(Z)\wih{\cdotM}\psi\\
&=&\wih{\nabla}_{Z}^M\psi-\frac{(-1)^{r+1} i}{2}Z\wih{\cdotM}\psi+\frac{(-1)^{r+1} i}{2}Z\wih{\cdotM}\psi-\frac{(-1)^{r+1}}{2}\wih{h}(Z)\wih{\cdotM}\psi_r+\frac{(-1)^{r+1}}{2}\wih{h}(Z)\wih{\cdotM}\psi_{r-1}\\
&=&\wih{\nabla}_{Z}^M\psi-\frac{(-1)^{r+1} i}{2}Z\wih{\cdotM}\psi+(-1)^{r+1} ip_-(Z)\wih{\cdotM}\psi_r+(-1)^{r+1} ip_+(Z)\wih{\cdotM}\psi_{r-1},
\ee
that is,
\[\left|\begin{array}{ll} 
\wih{\nabla}_Z\psi_r&=\big(\wih{\nabla}_{Z}^M\psi-\frac{(-1)^{r+1} i}{2}Z\wih{\cdotM}\psi\big)_r+(-1)^{r+1} ip_+(Z)\wih{\cdotM}\psi_{r-1}\\
\wih{\nabla}_Z\psi_{r-1}&=\big(\wih{\nabla}_{Z}^M\psi-\frac{(-1)^{r+1} i}{2}Z\wih{\cdotM}\psi\big)_{r-1}+(-1)^{r+1} ip_-(Z)\wih{\cdotM}\psi_r.\end{array}\right.\]
If one changes the Lorentzian metric $\wih{g}$ into $g:=-\wih{g}_{\wih{\xi}}\oplus \wih{g}_{\wih{\xi}^\perp}$, then one obtains a smooth Riemannian metric $g$ on $M$ and the triple $(M,g,\xi:=\wih{\xi})$ is a Riemannian flow with
\[ \left|\begin{array}{ll}\nabla_\xi^M\xi&=-\wih{\nabla}_{\wih{\xi}}^M\wih{\xi} \\ h&=-\wih{h}\\ \nabla&=\wih{\nabla}.\end{array}\right.\]
Moreover, the Clifford multiplications are related by 
\[\left|\begin{array}{ll}\xi\cdotM&=i\wih{\xi}\wih{\cdotM}\\ Z\cdotM&=Z\wih{\cdotM},\end{array}\right.\]
for all $Z\in\xi^\perp=\wih{\xi}^\perp$. Therefore the equations above become on $(M,g,\xi)$
\[\left|\begin{array}{ll}\nabla_{\xi}\psi_r&=\big(\wih{\nabla}_{\wih{\xi}}^M\psi-\frac{(-1)^{r+1} i}{2}\wih{\xi}\wih{\cdotM}\psi\big)_r-\frac{(-1)^{r} }{2}(n-2r)\xi\cdotM\psi_r\\ 
\nabla_{\xi}\psi_{r-1}&=\big(\wih{\nabla}_{\wih{\xi}}^M\psi-\frac{(-1)^{r+1} i}{2}\wih{\xi}\wih{\cdotM}\psi\big)_{r-1}+\frac{(-1)^{r} }{2}(n-2r)\xi\cdotM\psi_{r-1}\\
\nabla_Z\psi_r&=\big(\wih{\nabla}_{Z}^M\psi-\frac{(-1)^{r+1} i}{2}Z\wih{\cdotM}\psi\big)_r+(-1)^{r+1} ip_+(Z)\cdotM\psi_{r-1}\\
\nabla_Z\psi_{r-1}&=\big(\wih{\nabla}_{Z}^M\psi-\frac{(-1)^{r+1} i}{2}Z\wih{\cdotM}\psi\big)_{r-1}+(-1)^{r+1}ip_-(Z)\cdotM\psi_r.\end{array}\right.\]
Therefore, $\psi_r-i\epsilon\psi_{r-1}$ satisfies (\ref{eqspinKillinglorr}) if and only if $\psi$ is a $\frac{(-1)^{r+1}i}{2}$-Killing spinor on $(M,\wih{g},\wih{\xi})$.
\findemo

$ $\\

\noindent Round spheres provide examples of spin Sasakian manifolds where (\ref{eqspinKillingorig}) is fulfilled for the right $r$.

\begin{elemme}\label{lspinKillingsphere}
For any odd $n\geq 3$, the $(2n-1)$-dimensional round sphere $M$ with its canonical Sasakian and spin structures admits a $2\left(\begin{array}{c}n\\\frac{n+1}{2}\end{array}\right)$-dimensional space of sections of $\Sigma_{\frac{n+1}{2}}M\oplus\Sigma_{\frac{n-1}{2}}M\oplus\Sigma_{\frac{n-1}{2}}M\oplus\Sigma_{\frac{n-3}{2}}M$ satisfying {\rm (\ref{eqspinKillingorig})}. 
\end{elemme}

\noindent{\it Proof}:
Consider the standard embedding $\S^{2n-1}\subset\C^n$, with unit normal $\nu_x=x$ and hence Weingarten-endomorphism field $A=-\mathrm{Id}_{TM}$. Set $\xi:=-i\nu$. It is well-known that $(\S^{2n-1},g,\xi)$ is a Sasakian spin manifold with $h=-J$ on $\xi^\perp\subset TM$, where $J$ is the standard complex structure induced from $\C^n$.
Let $\psi\in\Sigma_r\C^n$ with $r\in\{0,1\ldots,n\}$ (i.e., $\wit{\Omega}\cdot\psi=i(2r-n)\psi$ where $\wit{\Omega}$ is the standard K\"ahler form of $\C^n$). If $r$ is even then $\psi\in\Sigma^+\C^n$. In that case the spinorial Gauss formula reads \[\nabla^M_X\varphi=\nabla_X^{\C^n}\varphi-\frac{1}{2}A(X)\cdotM\varphi\] so that the restriction of $\psi$ on $\S^{2n-1}$ satisfies $\nabla^M_X\psi=\frac{1}{2}X\cdotM\psi$, i.e., is a $\frac{1}{2}$-Killing spinor. If $r$ is odd, then $\psi\in\Sigma_-\C^n$. The spinorial Gauss formula for a section $\varphi\in\Sigma^-\C^n_{|_{\S^{2n-1}}}$, which can be identif\mbox{}ied with $\Sigma \S^{2n-1}$ provided we change the sign of the Clif\mbox{}ford multiplication, reads then
\[\nabla^M_X\varphi=\nabla_X^{\C^n}\varphi+\frac{1}{2}A(X)\cdotM\varphi\]
for every $X\in TM$. We deduce that $\nabla_X^M\psi=-\frac{1}{2}X\cdotM\psi$
for every $X\in TM$, that is, the restriction of $\psi$ to $\S^{2n-1}$ is a $-\frac{1}{2}$-Killing spinor.
To sum up, the restriction of a constant section $\psi\in\Sigma_r\C^n$ to $M:=\S^{2n-1}$ is a $\frac{(-1)^r}{2}$-Killing spinor on $M$.
Decompose such a $\psi$ into $\psi=\psi_r+\psi_{r-1}$, see (\ref{eqdecSr}). From Lemma~\ref{lequivspinKilling} and $\mathrm{rk}_{\C}(\Sigma_r\C^n)=\left(\begin{array}{c}n\\r\end{array}\right)$ we conclude.
\findemo
 $ $\\

\noindent The analog of $\S^{2n-1}$ in the Lorentzian context is the Anti-deSitter spacetime $\mathbb{H}^{2n-1}$, that can be defined by
\[ \mathbb{H}^{2n-1}:=\{z\in\C^n\,|\,-|z_0|^2+\sum_{j=1}^{n-1}|z_j|^2=-1\}.\]

\begin{elemme}\label{lspinKillingads}
For any odd $n\geq 3$, the $(2n-1)$-dimensional Anti-deSitter spacetime $M:=\mathbb{H}^{2n-1}$ with its induced Lorentzian Sasakian structure (with $\wih{\xi}_x=ix$ and $\wih{h}=J$) and induced spin structure admits an $\left(\begin{array}{c}n\\ r\end{array}\right)$-dimensional space of $\frac{(-1)^{r+1}i}{2}$-Killing spinors lying pointwise in $\Sigma_rM\oplus\Sigma_{r-1}M$.
In par\-ti\-cu\-lar, if one considers the (Riemannian) Sasakian metric given by $-\wih{g}_{\wih{\xi}}\oplus \wih{g}_{\wih{\xi}^\perp}$, where $\wih{g}$ is the canonical Lorentzian metric of sectional curvature $-1$, then $\mathbb{H}^{2n-1}$ admits a $2\left(\begin{array}{c}n\\\frac{n+1}{2}\end{array}\right)$-dimensional space of sections of $\Sigma_{\frac{n+1}{2}}M\oplus\Sigma_{\frac{n-1}{2}}M\oplus\Sigma_{\frac{n-1}{2}}M\oplus\Sigma_{\frac{n-3}{2}}M$ satisfying {\rm (\ref{eqspinKillinglor})}.
\end{elemme}

\noindent{\it Proof}: First recall that $M$ is a Lorentzian Sasakian manifold and simultaneously an $\S^1$-bundle with totally geodesic fibres over $\C\mathrm{H}^{n-1}$. Just as for the sphere, one can restrict spinors from $\C^n$ onto $M$ so that the following Gauss-Weingarten-formula holds for all $\psi\in C^\infty(\C^n,\Sigma_{2n})$ and all $X\in TM$:
\be 
\nabla_X^M\psi&=&-\frac{A(X)}{2}\cdot\nu\cdot\psi\\
&=&\left|\begin{array}{ll}\frac{iA(X)}{2}\cdotM\psi&\textrm{ if }\psi(x)\in\Sigma_{2n}^+\;\forall x\\ -\frac{iA(X)}{2}\cdotM\psi&\textrm{ if }\psi(x)\in\Sigma_{2n}^-\;\forall x,\end{array}\right.
\ee
where $A(X):=\wnabla_X\nu$ is the Weingarten endormorphism of $M$ in $\C^n$. Moreover, there still exists a $\wnabla$-parallel splitting $\Sigma_{2n}=\oplus_{r=0}^n\Sigma_{2n,r}$ where $\Sigma_{2n,r}:=\mathrm{Ker}(\wit{\Omega}\cdot-i(2r-n)\mathrm{Id})$ (with dimension $\left(\begin{array}{c}n\\ r\end{array}\right)$) and $\wit{\Omega}$ is the  K\"ahler form associated to the standard complex structure $J$ on $\wit{M}$. Choosing $\nu_x:=-x$ as unit normal on $M$, one has $A=-\mathrm{Id}_{TM}$, so that the restriction of any constant section of $\C^n\times\Sigma_{2n,r}$ onto $M$ provides a $\frac{(-1)^{r+1}i}{2}$-Killing spinor. Since again $\Sigma_r\wit{M}_{|_M}=\Sigma_rM\oplus\Sigma_{r-1}M$, the first statement follows. The second statement is a consequence of the first one together with Lemma~\ref{lequivspinKillinglor}. 
\findemo

$ $\\

\noindent The doubly warped product of Theorem~\ref{pclassifdwp}$.ii)$ corresponding to $M=\S^{2n-1}$ is the complement of a point in the complex hyperbolic space $\C\mathrm{H}^n$ with its canonical Fubini-Study metric of constant holomorphic sectional curvature $-4$ (compare with \cite[Satz 5.1]{BaierDiplA}).
Therefore we obtain a new description of the i\-ma\-gi\-na\-ry K\"ahlerian Killing spinors on $\C\mathrm{H}^n$ after the explicit one by K.-D. Kirchberg \cite[Sec. 3]{KirchbergAGAG93}.
Actually $\C\mathrm{H}^n$ is essentially the only example occurring in Theorem \ref{pclassifdwp}$.ii)$:

\begin{ethm}\label{pcarCHn}
For $n\geq 3$ odd let $(\wit{M}^{2n},\wit{g},\wit{J})$ be a K\"ahler doubly warped product as in {\rm Lemma \ref{lspindwp}} with $(M^{2n-1},\wih{g},\wih{\xi})$ complete, Sasakian, simply-connected, spin, $I=\R_+^\times$, $\rho=\sinh$ and $\sigma=\cosh$. Let $\wit{M}$ carry the induced spin structure and assume $(\wit{M}^{2n},\wit{g},\wit{J})$ admits a non-zero $i$-K\"ahlerian Killing spinor $(\psi,\phi)$.\\
Then $(\wit{M}^{2n},\wit{g},\wit{J})$ is holomorphically isometric to $\C\mathrm{H}^n\setminus\{x\}$ for some $x\in\C\mathrm{H}^n$.
\end{ethm}

\noindent{\it Proof}: It suffices to show that $(M^{2n-1},\wih{g},\wih{\xi})$ is $\S^{2n-1}$ with its standard Sasakian structure.
By assumption and Lemma \ref{lequivspinKilling}, the section $\varphi_{\frac{n+1}{2}}+\varphi_{\frac{n-1}{2}}$ is a $\frac{(-1)^{\frac{n+1}{2}}}{2}$-Killing spinor on $(M^{2n-1},\wih{g},\wih{\xi})$ lying pointwise in $\Sigma_{\frac{n+1}{2}}M\oplus\Sigma_{\frac{n-1}{2}}M$ and the section $\wit{\varphi}_{\frac{n-1}{2}}+\wit{\varphi}_{\frac{n-3}{2}}$ is a $-\frac{(-1)^{\frac{n+1}{2}}}{2}$-Killing spinor on $(M^{2n-1},\wih{g},\wih{\xi})$ lying pointwise in $\Sigma_{\frac{n-1}{2}}M\oplus\Sigma_{\frac{n-3}{2}}M$. At least one of them does not vanish. Now C. B\"ar's classification (see in particular \cite[Thm. 3]{Baer90}) implies that either $M=\S^{2n-1}$ or $M$ is a compact Einstein-Sasakian manifold with exactly one non-zero $\frac{1}{2}$- and one non-zero $-\frac{1}{2}$-Killing spinor. Moreover, each Killing spinor induces a parallel spinor on the Riemannian cone $\ovl{M}$ over $M$ \cite{Baer90}. But coming back to McK. Wang's classification of simply-connected complete Riemannian spin manifolds with parallel spinors, it turns out that, in the latter case, the reduced holonomy of $\ovl{M}$ is $\mathrm{SU}_n$ (where $n$ is its complex dimension) and the parallel spinors lie in $\Sigma_0\ovl{M}$ and $\Sigma_n\ovl{M}$ (see \cite[(i) p.61]{Wang89}), in particular not in $\Sigma_{\frac{n\pm1}{2}}\ovl{M}$. Thus only $\S^{2n-1}$ occurs.
\findemo

$ $\\

\noindent In case $M=\mathbb{H}^{2n-1}$ is equipped with its associated Riemannian Sasakian structure, the corresponding doubly warped product with $\rho=\cosh$ and $\sigma=\sinh$ has again constant holomorphic sectional curvature $-4$ by Lemma \ref{lholseccurv}.
It is actually the complement in $\C\mathrm{H}^n$ of some submanifold.
We conjecture that, up to covering, $\mathbb{H}^{2n-1}$ is the only Lorentzian Sasakian manifold having non-zero imaginary Killing spinors lying pointwise in the ``middle'' eigenspaces $\Sigma_{r}M$ (with $r\in\{\frac{n-3}{2},\ldots,\frac{n+1}{2}\}$) of the Clifford action of the transversal K\"ahler form.
If this happens, then only the complex hyperbolic space can occur as (simply-connected complete) example of doubly warped product in Theorem \ref{pclassifdwp}$.iii)$.
$ $\\

\section{Classification in a particular case}

\noindent In this section, we show that the structure of a doubly warped product can be recovered from the length function of a non-zero imaginary K\"ahlerian Killing spinor satisfying certain supplementary assumption on the K\"ahler manifold $\wit{M}$.
The following result can be seen as analogous to H. Baum's one \cite{Baum891} about imaginary Killing spinors of so-called \emph{type I}.
Recall for the next theorem that $V$ was defined by (\ref{eqdefV}).

\begin{ethm}\label{pclassiftype1}
Let $(\wit{M}^{2n},g,J)$ be a connected complete K\"ahler spin manifold carrying a non-zero $i$-K\"ahlerian Killing spinor $(\psi,\phi)$.
Assume $|\psi|=|\phi|$ and the existence of a real vector field $W$ on $\wit{M}$ together with a non-identically vanishing continuous function $\mu:\wit{M}\lra\C$ such that $W\cdot\psi=\mu\phi$.
Then the vector field $V$ has no zero, the K\"ahler manifold $(\wit{M}^{2n},g,J)$ is a doubly warped product as in {\rm Theorem~\ref{pclassifdwp}$.i)$} and $(\psi,\phi)$ comes from a transversally parallel spinor on $(M,\wih{g},\wih{\xi})$.
\end{ethm}

\noindent{\it Proof}: We construct a holomorphic isometry between $(\wit{M}^{2n},g,J)$ and some doubly warped product.
This isometry is provided by the flow of some vector field associated to the K\"ahlerian Killing spinor (compare with the case of imaginary Killing spinors \cite{Baum891}).\\
First note that, if $|\psi|=|\phi|$, then both $\psi$ and $\phi$ have no zero on $\wit{M}$. 
Because of $|W|\cdot |\psi|=|W\cdot\psi|=|\mu|\cdot|\phi|$, this already implies $|W|=|\mu|$ on $\wit{M}$.
Fix a neighbourhood $U$ of a point $x$ with $\mu(x)\neq 0$ for all $x\in U$.
It follows from the definition of $V$ that \begin{equation}\label{eq:muWV}\mu=2i\frac{g(p_+(W)_,V)}{|\phi|^2}\end{equation} 
on $U$, in particular $W(x)\neq 0$ and $V(x)\neq 0$ for all $x\in U$.
Now Cauchy-Schwarz inequality with $X=V$ in (\ref{eqdefV}) gives $|V|\leq|\psi|\cdot|\phi|$ on $\wit{M}$.
With (\ref{eq:muWV}) we deduce that
\be 
|\mu|^2&=&\frac{|V|^2}{|\phi|^4}\Big(g(W,\frac{V}{|V|})^2+g(W,\frac{J(V)}{|V|})^2\Big)\\
&\leq&\frac{|V|^2|W|^2}{|\phi|^4}\\
&\leq&|W|^2
\ee
on $U$, which together with $|\mu|=|W|$ provides $|V|=|\phi|^2$.
By the equality case in Cauchy-Schwarz inequality, we obtain $V\cdot\psi=i|V|\phi$ and $V\cdot\phi=i|V|\psi$ on $U$.
This identity holds on $\wit{M}$ because of the analyticity of all objects involved (by definition, $\psi$ is anti-holomorphic and $\phi$ is holomorphic).
This in turn implies $|V|=|\phi|^2$ on $\wit{M}$, in particular $\{V=0\}=\varnothing$ and $\frac{V}{|V|}\cdot\psi=i\phi$ as well as $\frac{V}{|V|}\cdot\phi=i\psi$ on $\wit{M}$.\\
Next we look at the level hypersurfaces $M_r:=\{x\in\wit{M},\,|\phi(x)|=r\}$ (with $r\in\R_+^\times$) which, if non-empty, are smooth because of $\{V=0\}=\varnothing$ and Proposition \ref{pKirch}.
A unit normal to $M_r$ is given by $\nu:=\frac{V}{|V|}$ and the associated Weingarten endomorphism field is 
\be 
A(X)&:=&-\wnabla_X\nu\\
&=&-\frac{1}{|V|}\left(\wnabla_X V-g(\wnabla_X V,\frac{V}{|V|})\frac{V}{|V|}\right)
\ee
for every $X\in\nu^\perp$.
Setting $\xi:=-J(\nu)$ (note that the vector field $\xi$ is pointwise tangent to $M_r$), using $\nu\cdot\psi=i\phi$ and Proposition \ref{pKirch}$.ii)$, we compute, for all $X,Y\in\nu^\perp$,
\be 
g(A(X),Y)&=&-\frac{1}{|V|}g(\wnabla_X V,Y)\\
&=&-\frac{1}{|V|}\Re e\pa{\langle p_-(X)\cdot\phi,p_-(Y)\cdot\phi\rangle+\langle p_+(X)\cdot\psi,p_+(Y)\cdot\psi\rangle}\\
&=&-\frac{1}{|V|}\Re e\pa{\langle p_-(X)\cdot\nu\cdot\psi,p_-(Y)\cdot\nu\cdot\psi\rangle+\langle p_+(X)\cdot\psi,p_+(Y)\cdot\psi\rangle}\\
&=&-\frac{1}{|V|}\Re e\Big(-\la p_-(X)\cdot\nu\cdot\psi,\nu\cdot p_-(Y)\cdot\psi\ra-2\ovl{g(\nu,p_-(Y))}\la p_-(X)\cdot\nu\cdot\psi,\psi\ra\\
& &\phantom{-\frac{1}{|V|}\Re e\Big(}+\langle p_+(X)\cdot\psi,p_+(Y)\cdot\psi\rangle\Big)\\
&=&-\frac{1}{|V|}\Re e\Big(\la\nu\cdot p_-(X)\cdot\psi,\nu\cdot p_-(Y)\cdot\psi\ra+2g(\nu,p_-(X))\la\psi,\nu\cdot p_-(Y)\cdot\psi\ra\\
& &\phantom{-\frac{1}{|V|}\Re e\Big(}-2\ovl{g(\nu,p_-(Y))}\la p_-(X)\cdot\nu\cdot\psi,\psi\ra+\langle p_+(X)\cdot\psi,p_+(Y)\cdot\psi\rangle\Big)\\
&=&-\frac{1}{|V|}\Re e\Big(\la X\cdot\psi,Y\cdot\psi\ra+ig(\nu,J(X))\la\psi,\nu\cdot p_-(Y)\cdot\psi\ra+ig(\nu,J(Y))\la p_-(X)\cdot\nu\cdot\psi,\psi\ra\Big)\\
&=&-\frac{1}{|V|}\Big(|\psi|^2g(X,Y)+g(\nu,J(X))\Re e(\underbrace{\la\phi,p_-(Y)\cdot\psi\ra}_0)-g(\nu,J(Y))\Re e(\la p_-(X)\cdot\phi,\psi\ra)\Big)\\
&=&-\frac{1}{|V|}\Big(|\psi|^2g(X,Y)+g(\nu,J(Y))g(J(X),V)\Big)\\
&=&-(g(X,Y)+g(\xi,X)g(\xi,Y)),
\ee
that is, $A=-\mathrm{Id}_{TM_r}-\xi^\flat\otimes\xi$.
In particular, the Gau\ss{}-Weingarten formula for the inclusion $M_r\subset\wit{M}$ reads $\wnabla_XY=\nabla_X^{M_r}Y-(g(X,Y)+g(\xi,X)g(\xi,Y))\nu$ for all vector fields $X,Y$ tangent to $M_r$.\\
We begin with the reconstruction of the doubly warped product structure of Theorem \ref{pclassifdwp}$.i)$.
From $A(\xi)=-2\xi$, we deduce that $A(J(V))=-2J(V)$, hence $\wnabla_{J(V)}\nu=2J(V)$.
Proposition \ref{pKirch}$.ii)$ gives
\[J(V)(|V|)=\frac{g(\wnabla_VV,J(V))}{|V|}=\frac{1}{|V|}\Re e\pa{\langle p_-(V)\cdot\phi,p_-(J(V))\cdot\phi\rangle+\langle p_+(V)\cdot\psi,p_+(J(V))\cdot\psi\rangle}=0.\]
Therefore $\wnabla_{J(V)}V=2|V|J(V)$, that is, $\wnabla_VV=2|V|V$ using $\wnabla_{J(X)}V=J(\wnabla_XV)$ for all $X$.
This implies for the commutator of $\xi$ and $\nu$ (which we need later for the identification of the metric and of the Sasakian structure)
\begin{eqnarray}\label{eqcrochetxinu}
\nonumber[\xi,\nu]&=&-[J(\nu),\nu]\\
\nonumber&=&-[\frac{J(V)}{|V|},\frac{V}{|V|}]\\
\nonumber&=&-\frac{1}{|V|}\underbrace{J(V)(\frac{1}{|V|})}_0 V-\frac{1}{|V|}[\frac{J(V)}{|V|},V]\\
\nonumber&=&\frac{1}{|V|}V(\frac{1}{|V|})J(V)-\frac{1}{|V|^2}[J(V),V]\\
\nonumber&=&-\frac{g(\wnabla_VV,V)}{|V|^3}J(\frac{V}{|V|})-\frac{1}{|V|^2}J(\underbrace{[V,V]}_0)\\
&=&2\xi.
\end{eqnarray}
We show now that each (non-empty) $(M_r,g_{|_{M_r}},\xi_{|_{M_r}})$ is Sasakian.
For every $X\in TM_r$, one has
\be 
\wnabla_X\xi&=&-\wnabla_X(J(\nu))\\
&=&-J(\wnabla_X\nu)\\
&=&J(A(X))\\
&=&-J(X)-g(\xi,X)\nu,
\ee
so that $\wnabla_\xi\xi=-2\nu$, from which $\nabla_\xi^{M_r}\xi=0$ follows and, for every $Z\in\{\xi,\nu\}^\perp$, the identity $\wnabla_Z\xi=-J(Z)$ implies $\nabla_Z^{M_r}\xi=-J(Z)$.
In particular, $\xi_{|_{M_r}}$ defines a minimal Riemannian flow on $(M_r,g_{|_{M_r}})$ and $h=-J$ is an almost Hermitian structure on $\xi^\perp\subset TM_r$.
It remains to show that $h$ - or, equivalently, $J$ - is transversally parallel on $\xi^\perp$.
Recall that, from the definition of the transversal covariant derivative $\nabla$ one has, for all sections $Z,Z'$ of $\xi^\perp$,
\be 
\nabla_\xi Z&=&\nabla_\xi^{M_r} Z-h(Z)\\
&=&\wnabla_\xi Z-g(A(\xi),Z)\nu+J(Z)\\
&=&\wnabla_\xi Z+J(Z)
\ee
and
\be 
\nabla_Z Z'&=&\nabla_Z^{M_r} Z'+g(h(Z),Z')\xi\\
&=&\wnabla_Z Z'-g(A(Z),Z')\nu-g(J(Z),Z')\xi\\
&=&\wnabla_Z Z'+g(Z,Z')\nu-g(J(Z),Z')\xi,
\ee
from which one deduces that
\be 
(\nabla_\xi J)(Z)&=&\nabla_\xi (J(Z))-J(\nabla_\xi Z)\\
&=&\wnabla_\xi (J(Z))-Z-J(\wnabla_\xi Z)+Z\\
&=&0
\ee
and
\be 
(\nabla_Z J)(Z')&=&\nabla_Z(J(Z'))-J(\nabla_Z Z')\\
&=&\wnabla_Z (J(Z'))+g(Z,J(Z'))\nu-g(J(Z),J(Z'))\xi\\
& &-J(\wnabla_Z Z')+g(Z,Z')\xi+g(J(Z),Z')\nu\\
&=&0,
\ee
i.e., $\nabla J=0$, which proves that $(M_r,g_{|_{M_r}},\xi_{|_{M_r}})$ is Sasakian.\\
We come to the holomorphic isometry.
Denote $M:=M_1$, $\wih{g}:=g_{|_M}$ and $\wih{\xi}:=\xi_{|_M}$.
Up to rescaling $(\psi,\phi)$ by a positive constant (this does not influence both conditions on $(\psi,\phi)$), we may assume that $M\neq\varnothing$.
Let $F_t^\nu$ be the flow of $\nu$ on $\wit{M}$.
The vector field $\nu$ is complete since $\nu$ is bounded and $(\wit{M},g)$ is complete.
Consider the map
\be 
F:M\times \R&\lra&\wit{M}\\
(x,t)&\lmt&F_t^\nu(x).
\ee
We first show that $F$ is a diffeomorphism.
If $F_t^\nu(x)=F_{t'}^\nu(x')$ for some $t,t'\in\R$ and $x,x'\in M$, then $x$ and $x'$ lie on the same integral curve of $\nu$.
Let now $c$ be any integral curve of $\nu$ on $\wit{M}$ with $c(0)\in M$ and set $f(t):=|V|_{c(t)}$ (note that $f$ {\sl a priori} depends on the curve and in particular on the chosen starting point).
Then $f$ is smooth with first derivative given by $f'(t)=\frac{g(\wnabla_VV,V)}{|V|^2}(c(t))=2|V|_{c(t)}=2f(t)$ for all $t$, so that $f=f(0)e^{2t}=e^{2t}$.
This has several consequences.
On the one hand, $f$ is injective, so that $c$ meets $M$ at most once, hence $x=x'$ and $t=t'$, which proves the injectivity of $F$.
On the other hand, $f$ does {\sl a posteriori} not depend on the chosen starting point on $M$, in particular $F_t^\nu$ preserves the foliation by the level hypersurfaces $M_r$ of $|\phi|$ and hence the orthogonal splitting $TM_r\oplus\R\nu$.
Together with the surjectivity of $f:\R\rightarrow\R_+^\times$, we obtain that of $F$ and the pointwise invertibility of the differential of $F$.
Therefore $F$ is a diffeomorphism.\\
Next we determine the metric $F^*g$.
The map $F$ sends $\frac{\partial}{\partial t}$ onto $\nu$,
so that obviously $F^*g(\frac{\partial}{\partial t},\frac{\partial}{\partial t})=1$.
The preceding considerations also yield $F^*g(\frac{\partial}{\partial t},X)=0$ for all $t\in\R$ and $X\in TM$.
Since 
\[\frac{\partial}{\partial s}(F_s^\nu)_*\xi_{|_{s=t}}=(F_t^\nu)_*\frac{\partial}{\partial s}(F_s^\nu)_*\xi_{|_{s=0}}=(F_t^\nu)_*[\xi,\nu]\bui{=}{(\ref{eqcrochetxinu})}2(F_t^\nu)_*\xi,\]
we have
\begin{equation}\label{eqFtxi}
(F_t^\nu)_*\xi=e^{2t}\xi 
\end{equation}
for every $t\in \R$. Moreover, the Lie derivative of $g$ in direction of $\nu$ is given for all $X,Y\in\nu^\perp$ by
\be 
(\mathcal{L}_\nu g)(X,Y)&=&g(\wnabla_X\nu,Y)+g(\wnabla_Y\nu,X)\\
&=&-2g(A(X),Y)\\
&=&2(g(X,Y)+g(\xi,X)g(\xi,Y)),
\ee
that is, $(\mathcal{L}_\nu g)_{|_{\nu^\perp}}=2(g+\xi^\flat\otimes\xi^\flat)$.
The identity $\frac{\partial}{\partial s}(F_s^\nu)^*g_{|_{s=t}}=(F_t^\nu)^*\mathcal{L}_\nu g $ provides, for any $X,Y\in TM$ and $t\in\R$
\begin{eqnarray}\label{eqdFtg} 
\nonumber\frac{\partial}{\partial s}((F_s^\nu)^*g(X,Y))_{|_{s=t}}&=&(\frac{\partial}{\partial s}(F_s^\nu)^*g{}_{|_{s=t}})(X,Y)\\
\nonumber&=&\{(F_t^\nu)^*\mathcal{L}_\nu g\}(X,Y)\\
\nonumber&=&\mathcal{L}_\nu g((F_t^\nu)_*X,(F_t^\nu)_*Y)\circ F_t^\nu\\
\nonumber&=&2\Big(g((F_t^\nu)_*X,(F_t^\nu)_*Y)+g(\xi,(F_t^\nu)_*X)g(\xi,(F_t^\nu)_*Y)\Big)\circ F_t^\nu\\
\nonumber&=&2\Big((F_t^\nu)^*g(X,Y)+(F_t^\nu)^*g((F_{-t}^\nu)_*\xi,X)(F_t^\nu)^*g((F_{-t}^\nu)_*\xi,Y)\Big)\\
&\bui{=}{\rm (\ref{eqFtxi})}&2\Big((F_t^\nu)^*g(X,Y)+e^{-4t}(F_t^\nu)^*g(\xi,X)(F_t^\nu)^*g(\xi,Y)\Big).
\end{eqnarray}
Since $(F_t^\nu)^*g(\xi,\xi)=g((F_t^\nu)_*\xi,(F_t^\nu)_*\xi)\circ F_t^\nu\bui{=}{\rm (\ref{eqFtxi})}(e^{4t}g(\xi,\xi))\circ F_t^\nu=e^{4t}$, we deduce from (\ref{eqdFtg}) that, for $X=\xi$,
\[\frac{\partial}{\partial s}((F_s^\nu)^*g(\xi,Y))_{|_{s=t}}=4(F_t^\nu)^*g(\xi,Y),\]
from which $(F_t^\nu)^*g(\xi,Y)=e^{4t}g(\xi,Y)$ follows.
In particular, $(F_t^\nu)^*g(\xi,Y)=0$ for every $Y\in\{\xi,\nu\}^\perp$.
For $X,Y\in\{\xi,\nu\}^\perp$, the identity (\ref{eqdFtg}) becomes
\[ \frac{\partial}{\partial s}((F_s^\nu)^*g(X,Y))_{|_{s=t}}=2(F_t^\nu)^*g(X,Y),\]
which implies $(F_t^\nu)^*g(X,Y)=e^{2t}g(X,Y)$.
To sum up, the pull-back metric on $M\times\R$ is given by
\[ F^*g=e^{2t}(e^{2t}\wih{g}_{\wih{\xi}}\oplus \wih{g}_{\wih{\xi}^\perp})\oplus dt^2,\]
where $\wih{g}_{\wih{\xi}}=\wih{\xi}^\flat\otimes\wih{\xi}^\flat=\wih{g}(\wih{\xi},\cdot)\otimes\wih{g}(\wih{\xi},\cdot)$ and, as in the beginning of this section, $\wih{g}_{\wih{\xi}^\perp}$ denotes the restriction of $\wih{g}$ onto the subspace $\{\wih{\xi},\frac{\partial}{\partial t}\}^\perp\subset TM$. 
Hence the map $F$ provides an isometry with the doubly warped product of Theorem \ref{pclassifdwp}$.i)$.
This isometry pulls the spin structure of $\wit{M}$ back onto the product spin structure of $M\times\R$, where $M$ carries the spin structure induced by its embedding in $\wit{M}$.
It remains to show that $F$ identifies the complex structures.
This follows from the definition of the complex structure on the doubly warped product $M\times\R$ (see Lemma \ref{lstructKaehl}), from $(F_t^\nu)_*\nu=\nu$, $(F_t^\nu)_*(e^{-2t}\wih{\xi})=\xi$ and from $[J(Z),\nu]=\wnabla_{J(Z)}\nu-\wnabla_\nu J(Z)=-A(J(Z))-J(\wnabla_\nu Z)=J(Z)-J(\wnabla_\nu Z)=J([Z,\nu])$ for every section $Z$ of $\{\xi,\nu\}^\perp$ (use the computation of $A$ above).\\
{\sl Last but not the least}, the identity $\nu\cdot\psi=i\phi$ implies that $\phi$ (or, equivalently, $\psi$) is transversally parallel on $(M,\wih{g},\wih{\xi})$ by Theorem \ref{pclassifdwp}$.i)$.
This concludes the proof of Theorem \ref{pclassiftype1}.
\findemo

$ $\\

\noindent It is important to note that only the condition $W\cdot\psi=\mu\phi$ for some real vector field $W$ is restrictive, since by \cite[Thm. 11]{KirchbergAGAG93} the identity $|\psi|=|\phi|$ can always be assumed.\\

\noindent We conjecture that the examples of Section \ref{s:dwp} describe all K\"ahler spin manifolds admitting non-trivial imaginary K\"ahlerian Killing spinors.
This will be the object of a forthcoming paper.

\providecommand{\bysame}{\leavevmode\hbox to3em{\hrulefill}\thinspace}


\begin{thebibliography}{10}
\addcontentsline{toc}{chapter}{References}

\bibitem{BaierDiplA}
P.D. Baier, \emph{\"Uber den Diracoperator auf Mannigfaltigkeiten mit Zylinderenden}, Diplomarbeit, Universit\"at Freiburg, 1997.

\bibitem{Baer90}
C. B\"ar, \emph{Real Killing spinors and holonomy}, Comm. Math. Phys. \textbf{154} (1993), no. 3, 509--521.

\bibitem{Baum891}
H. Baum, \emph{Complete Riemannian manifolds with imaginary Killing spinors}, Ann. Glob. Anal. Geom. \textbf{7} (1989), 205--226.

\bibitem{Baum892}
\bysame, \emph{Odd-dimensional Riemannian manifolds admitting imaginary Killing spinors},  Ann. Glob. Anal. Geom. \textbf{7} (1989), 141--153.



%
%

\bibitem{Beauville83}
A. Beauville,  \emph{Vari\'et\'es K\"ahleriennes dont la premi\`ere classe de Chern est nulle}, J. Diff. Geom. \textbf{18} (1983); no. 4, 755--782 (1984).

\bibitem{GinHabibspKTgeom}
N. Ginoux and G. Habib, \emph{Geometric aspects of transversal Killing spinors on Riemannian flows}, Abh. Math. Sem. Univ. Hamburg \textbf{78} (2008), 69--90.


\bibitem{Habibthese}
G. Habib, \emph{Tenseur d'impulsion-\'energie et feuilletages}, PhD thesis, Institut \'Elie Cartan - Universit\'e Henri Poincar\'e, Nancy (2006).

\bibitem{HabibEM}
\bysame, \emph{Energy-Momentum tensor on foliations}, J. Geom. Phys. \textbf{57} (2007), no. 11, 2234--2248.


%

\bibitem{KirchbergAGAG93}
K.-D. Kirchberg, \emph{Killing spinors on K\"ahler manifolds}, Ann. Glob. Anal. Geom. \textbf{11} (1993), 141--164.

\bibitem{Kocomplexgeom}
S. Kobayashi, \emph{Transformation groups in differential geometry}, Ergebnisse der Mathematik und ihrer Grenzgebiete \textbf{70} (1972), Springer-Verlag.

\bibitem{KoNo}
S. Kobayashi, K. Nomizu, \emph{Foundations of differential geometry}, Vol. I-II, Wiley-Interscience Publication, New-York, 1963-1969.


%

\bibitem{Moroi95}
A. Moroianu, \emph{La premi\`ere valeur propre de l'op\'erateur de Dirac sur les vari\'et\'es k\"ahleriennes compactes},
Commun. Math. Phys. \textbf{169} (1995), 373--384.

\bibitem{MoroiHodge}
\bysame, \emph{Spineurs et vari\'et\'es de Hodge}, Rev. Roumaine Math. Pures Appl. \textbf{43} (1998), no. 5-6, 615--626.

%

\bibitem{Wang89}
McK. Wang, \emph{Parallel spinors and parallel forms}, Ann. Glob. Anal. Geom. \textbf{7} (1989), 59--68.
\end{thebibliography}
\end{document}